\documentclass[journal]{IEEEtran}

\usepackage{xspace,amssymb,epsfig,subfigure,syntonly,bm,bbm}
\usepackage{epsfig,amsmath,color}
\usepackage{xspace,syntonly}
\usepackage{mathtools}

\usepackage{algorithm,algorithmic}

\newcommand{\utwi}[1]{\mbox{\boldmath $#1$}}

\newcommand{\Prob}{{\textrm{Pr}}}

\newcommand{\cD}{{\cal D}}

\newcommand{\cL}{{\cal{L}}}
\newcommand{\cN}{{\cal N}}
\newcommand{\cP}{{\cal P}}

\newcommand{\cA}{{\cal A}}

\newcommand{\cT}{{\cal T}}

\newcommand{\cE}{{\cal E}}
\newcommand{\cI}{{\cal I}}
\newcommand{\cF}{{\cal F}}

\newcommand{\cU}{{\cal U}}

\newcommand{\cH}{{\cal H}}
\newcommand{\cV}{{\cal V}}

\newcommand{\cX}{{\cal X}}

\newcommand{\bd}{{\bf d}}
\newcommand{\be}{{\bf e}}

\newcommand{\bff}{{\bf f}}

\newcommand{\bs}{{\bf s}}
\newcommand{\bx}{{\bf x}}

\newcommand{\bv}{{\bf v}}
\newcommand{\bi}{{\bf i}}
\newcommand{\bz}{{\bf z}}

\newcommand{\bA}{{\bf A}}

\newcommand{\bI}{{\bf I}}

\newcommand{\bZ}{{\bf Z}}

\newcommand{\bY}{{\bf Y}}
\newcommand{\bV}{{\bf V}}

\newcommand{\bepsilon}{{\utwi{\epsilon}}}
\newcommand{\bgamma}{{\utwi{\gamma}}}

\newcommand{\bnu}{{\utwi{\nu}}}

\newcommand{\bPhi}{{\utwi{\Phi}}}
\newcommand{\bxi}{{\utwi{\xi}}}

\newcommand{\bvarphi}{{\utwi{\varphi}}}
\newcommand{\bchi}{{\utwi{\chi}}}
\newcommand{\bmu}{{\utwi{\mu}}}

\newcommand{\biota}{{\utwi{\iota}}}
\newcommand{\bsigma}{{\utwi{\sigma}}}
\newcommand{\bell}{{\utwi \ell}}

\newcommand{\bmA}{\bm{A}}

\newcommand{\bbmM}{\bar{\bm{M}}}

\DeclarePairedDelimiter{\ceil}{\lceil}{\rceil}

\begin{document}

\newtheorem{definition}{Definition}
\newtheorem{remark}{Remark}
\newtheorem{proposition}{Proposition}
\newtheorem{lemma}{Lemma}

%------------------------------------------------------------------------------
% Title.
%------------------------------------------------------------------------------

\title{Risk-Constrained Microgrid Reconfiguration \\
Using Group Sparsity}

\author{Emiliano Dall'Anese, \emph{Member, IEEE}, and Georgios B. Giannakis, \emph{Fellow, IEEE}
\thanks{\protect\rule{0pt}{0.5cm}%
Paper submitted in June 6, 2013; revised on April 9, 2014; accepted May 16, 2014. 
This work was supported by the Inst. of Renewable Energy and the
Environment (IREE) grant no. RL-0010-13, Univ. of Minnesota, and by the 
grants  NSF ECCS 1002180, NSF ECCS 1202135, NSF AST 1247885, and NSF EARS 1343248.
Authors are with the Digital Technology Center and the Dept. of ECE, University of Minnesota, 200 Union Street SE, Minneapolis, MN 55455, USA. 
E-mails: {\tt \{emiliano, georgios\}@umn.edu} 
}
}

% The paper headers
\markboth{IEEE TRANSACTIONS ON SUSTAINABLE ENERGY}%
{Dall'Anese and Giannakis: Risk-Constrained Microgrid Reconfiguration Using Group Sparsity }

% make the title area
\maketitle

\maketitle

%%%%%%%%%%%%%%%%%%%%%%%%%%%%%%%%%%%%%%%%%%%%%%%
% Abstract
%%%%%%%%%%%%%%%%%%%%%%%%%%%%%%%%%%%%%%%%%%%%%%%
\begin{abstract}
The system reconfiguration task is considered for existing power distribution systems and microgrids, in the presence of  renewable-based generation and load foresting errors. The system topology is obtained by solving a chance-constrained optimization problem, where loss-of-load (LOL) constraints and Ampacity limits of the distribution lines are enforced. Similar to various distribution system reconfiguration renditions, solving the resultant problem is computationally prohibitive due to the presence of binary line selection variables. Further, lack of closed form expressions for the joint probability distribution of forecasting errors hinders tractability of LOL constraints. Nevertheless, a \emph{convex} problem re-formulation is developed here by resorting to a scenario approximation technique, and by leveraging the underlying group-sparsity attribute of currents flowing on distribution lines equipped with tie and sectionalizing switches. The novel convex LOL-constrained reconfiguration scheme can also afford a distributed solution using the alternating direction method of multipliers, to address the case where multi-facilities are managed autonomously from the rest of the system. 
\end{abstract}


%%%%%%%%%%%%%%%%%%%%%%%%%%%%%%%%%%%%%%%%%%%%%%
\section{Introduction}
\label{sec:Introduction}
%%%%%%%%%%%%%%%%%%%%%%%%%%%%%%%%%%%%%%%%%%%%%%

Distributed energy resources (DERs) are critical modules of existing power distribution systems and future microgrids, 
and one of the driving forces toward transforming today's distribution grid into a sustainable, scalable, and efficient one~\cite{Hatziargyriou-PESMag}. DERs include small-scale controllable power sources such as diesel
generators and micro combined heat and power (microCHP) units, as well as renewable energy sources (RESs), with photovoltaic (PV) systems and small-scale wind turbines as prime examples. DERs bring generation closer to the end user, offer environment-friendly advantages, and can also provide ancillary services~\cite{Hatziargyriou-PESMag}. 

RES generation is stochastic, non-dispatchable, and challenging to predict accurately in real-time~\cite{Bacher2009,Lorenz09,Tsikalakis06}. Although numerical weather forecasts yield reasonably reliable predictions of the average solar irradiance and wind speed over intervals of say 10-15 minutes~\cite{Bacher2009,Tsikalakis06}, the instantaneous power available may unexpectedly fluctuate around its forecasted value due to e.g., variable cloud coverage and gusts of wind. An additional potential source of uncertainty is load forecasting errors~\cite{Hodge2013}, especially in the presence of stochastic elastic load demand patterns, such as those corresponding to electric vehicles. In fact, customers may decide to start charging vehicles at their convenience, rather than relying on aggregator policies. These sources of uncertainty in RES generation and load demand may lead microgrids to operate possibly far from the expected regime, where steady-state variables are fine-tuned based on load, solar, and wind predictions~\cite{Khodr09,Dallanese-TSG13}. Potential consequences include, for instance, loss of load (LOL) at one or more nodes, and line overheating which, in turn, may trigger outages. Thus, for both short- and long-term microgrid operation planning (from a few minutes to hours ahead), it is essential to account for uncertain RES generation and load profiles, in order to ensure a reliable power delivery microgrid-wide, make risk-limiting operational decisions, and facilitate the penetration of RESs in large-scale~\cite{Hatziargyriou-PESMag}. 

The impact of intermittent RES generation on the economic dispatch task was considered in~\cite{Liu2010,Villanueva12,Zhang-isgt13} (and pertinent references therein).  A rolling horizon strategy for energy management in microgrids with renewables was proposed in~\cite{Palma-Behnke13}, where forecasted ambient conditions were also accounted for.
However, conventional economic dispatch strategies are oblivious to electrical network constraints and power losses, which may play a critical role in determining the supply-demand (im)balance. The effects of uncertain generation on the electrical network were assessed in~\cite{Morales2010} and~\cite{Ahmadi2011} using a probabilistic power flow approach to test the system functionalities over a variety of operational conditions; see also~\cite{Cortes-Carmona10}, where fuzzy arithmetic was adopted for a probabilistic DC load flow problem. Chance-constrained optimal power flow (OPF) formulations (using a DC approximation) were considered in~\cite{Chertkov13}. In the distribution system reconfiguration context, a probabilistic load flow scheme was employed in~\cite{Celli02} to identify the distribution network configuration that is more likely to adhere to thermal limits.  

The microgrid reconfiguration problem under uncertain load and RES generation is considered in this paper. The novel approach seeks the configuration that is optimal according to a well defined criterion, while ensuring \emph{satisfaction of the load} with arbitrarily \emph{high probability}, and strict adherence to maximum current limits on the distribution lines. Similar to past works on system reconfiguration without uncertainties~\cite{Khodr09,Baran89,Schmidt95,Cicoria04}, the formulated problem is hard to solve optimally and efficiently due to binary line selection variables and the nonlinear power flow relations. What is more, LOL constraints are intractable because no closed form expressions are typically available for the joint probability density function of the power supplied by multiple wind farms and PV systems. Nevertheless, a computationally affordable re-formulation can be derived by resorting to a Monte Carlo based scenario approximation technique~\cite{Calafiore,Zhang-isgt13}, and by exploiting the underlying group sparsity of currents and powers flowing over the
conductors of distribution lines equipped with switches. This group sparsity attribute enables re-casting the reconfiguration task using a constrained multidimensional shrinkage and thresholding operator (MSTO)~\cite{YuLi06,Wiesel11}. The upshot of the proposed approach is that the resultant formulation is \emph{convex} and \emph{sample-size-free}. Unlike competing alternatives that require solving a nonconvex power flow problem per sample~\cite{Morales2010,Ahmadi2011,Celli02}, the proposed approach entails solving a constrained MSTO problem with a single supply-demand balance constraint per phase and node. 

To accommodate microgrids which include single- or multi-facilities that are managed independently from the rest of the network, the proposed reconfiguration is solved in a decentralized fashion by resorting to the so-called alternating direction method of multipliers (ADMM)~\cite[Sec.~3.4]{BeT89}. In the power systems context, ADMM was employed in~\cite{Kekatos13} to estimate the state of transmission systems distributedly, in~\cite{Boyd_messPass,Erseghe14,Dallanese-TEC14} to derive distributed OPF solvers for balanced power systems, and in~\cite{Dallanese-TSG13} to solve the OPF problem for unbalanced distribution systems in a decentralized fashion. Here, the approach is tailored for the microgrid reconfiguration problem. 
The novel decentralized reconfiguration algorithm entails a two-way communication between the microgrid manager and the area controllers.

To summarize, the main contributions of this work are as follows:  \emph{i)} a novel distribution system reconfiguration problem in the presence of load and renewable generation forecasting errors is formulated, where LOL probability constraints and maximum current limits are enforced; \emph{ii)} a computationally affordable convex relaxation is derived by resorting to a Monte Carlo based scenario approximation technique~\cite{Calafiore}, along with the sparsity-promoting regularization techniques first utilized in~\cite{Dallanese-TPD13} to bypass line selection variables; and, \emph{iii)} an ADMM-based  algorithm is developed to solve the proposed risk-aware system reconfiguration problem in a decentralized fashion.\footnote{{\it Notation:} Upper (lower) boldface
letters will be used for matrices (column vectors); $(\cdot)^\cT$ for transposition; $(\cdot)^*$ for complex-conjugate; and, $(\cdot)^\cH$ for complex-conjugate transposition; $\Re\{\cdot\}$ will denote the real part; $\Im\{\cdot\}$ the imaginary part; and, $j := \sqrt{-1}$ the imaginary unit; $|\cP|$ the cardinality of set $\cP$; and, $\mathbb{R}^{N}$ and $\mathbb{C}^{N}$ the space of the $N \times 1$ real and complex vectors, respectively. Given vector $\bv$ and matrix $\bV$, $[\bv]_{\cP}$ will denote a $|\cP| \times 1$ sub-vector containing entries of $\bv$ indexed by the set $\cP$, and $[\bV]_{\cP_1,\cP_2}$ the $|\cP_1| \times |\cP_2|$ sub-matrix with row and column indexes described by $\cP_1$ and $\cP_2$. Further, $\|\bv \|_2 := \sqrt{\bv^{\cH}\bv}$ will stand for the $\ell_2$-norm of $\bv$; and $\mathbf{0}_{M\times N}$, $\mathbf{1}_{M\times N}$ for $M \times N$ matrices with all zeroes and ones, respectively; and, $\ceil{a}$ for the
smallest integer greater than or equal to $a$. Finally, $\preceq$ and $\succeq$ are element-wise inequalities, and $\Prob\{A\}$ will denote the probability of event $A$.}
%

%%%%%%%%%%%%%%%%%%%%%%%%%%%%%%%%%%%%%%%%%%%%%%
\section{Preliminaries and problem statement}
\label{sec:Formulation}
%%%%%%%%%%%%%%%%%%%%%%%%%%%%%%%%%%%%%%%%%%%%%%

Consider modeling a microgrid as a directed graph $(\cN, \cE)$, where $N$ nodes are collected in the set $\cN := \{1,\ldots,N\}$, and overhead or underground lines are represented by the set of directed edges $\cE := \{(m,n)\} \subset \cN  \times \cN$.
Let node $1$ represent the point of common coupling (PCC). Define as $\cP_{mn} \subseteq \{a, b, c\}$ and $\cP_{n} \subseteq \{a, b, c\}$ the phases of line $(m,n) \in \cE$ and node $n \in \cN$, respectively. Let $V_n^{\phi} \in \mathbb{C}$ be the phasor representation of the complex line-to-ground voltage at node $n \in \cN$ of phase $\phi \in \cP_n$, and likewise $I_n^{\phi} \in \mathbb{C}$ for the current injected. Lines $(m,n) \in \cE$ are modeled as $\pi$-equivalent components~\cite[Ch.~6]{Kerstingbook}, with $\bZ_{mn} \in \mathbb{C}^{|\cP_{mn}| \times |\cP_{mn}|}$ and $\bY_{mn} \in \mathbb{C}^{|\cP_{mn}| \times |\cP_{mn}|}$ denoting the  phase impedance and shunt admittances matrices, respectively. Typically, $\bZ_{mn}$ is symmetric (but not Hermitian), full rank, and has non-zero diagonal elements~\cite{testfeeder}. On the other hand, $\bY_{mn}$ is typically diagonal~\cite{Kerstingbook,testfeeder}. Notice that values of the diagonal entries of $\bY_{mn}$ are on the order of $10-100$ micro Siemens per mile~\cite{testfeeder}; thus, similar to existing works on distribution system reconfiguration~\cite{Khodr09,Celli02,Baran89,Schmidt95,Bolognani13} (and references therein), the effects of shunt admittance matrices are neglected in this paper. However, one way to account for possible non-negligible effects of  shunt admittances will be described in Section~\ref{sec:dsr_sparsity}.

Changes in the microgrid topology are effected by opening or closing tie and sectionalizing line switches. Thus, collect in the subset $\cE_{R} \subset \cE$ the lines equipped with controllable switches, and let the binary variable $x_{mn} \in \{0,1\}$ indicate whether line $(m,n) \in \cE_R$ is used ($x_{mn} = 1$) or not ($x_{mn} = 0$); see, e.g.~\cite{Khodr09,Baran89,Schmidt95}. Set $\cE_R$ clearly includes the branch connecting the microgrid to the PCC; the microgrid  operates in a grid-connected mode when this branch is used, and in an islanded setup otherwise~\cite{Hatziargyriou-PESMag}. With this notation, the $|\cP_{mn}| \times 1$ vector $\bi_{mn} := [\{I_{mn}^\phi, \phi \in \cP_{mn} \}]^\cT$ collecting the currents flowing on line $(m,n) \in \cE_R$ can be collectively expressed as  
\begin{align}
\label{line_selection}
\bi_{mn} =   x_{mn} \bZ_{mn}^{-1}  \left( [\bv_m]_{\cP_{mn}} - [\bv_n]_{\cP_{mn}} \right), \,\, x_{mn} \in \{0,1\}
\end{align}
where $\bv_{m} := [\{V_{m}^{\phi}, \phi \in \cP_{m}\}]^\cT$. Clearly, the counterpart of~\eqref{line_selection} for lines $(m,n) \in \cE \backslash \cE_R$ reads $\bi_{mn} =  \bZ_{mn}^{-1}  \left( [\bv_m]_{\cP_{mn}} - [\bv_n]_{\cP_{mn}} \right)$. Line currents $\{I_{mn}^\phi \}$ and injected currents $\{I_{n}^\phi\}$ abide by Kirchhoff's current law, which can be written per phase $\phi$ and node $n$ as
\begin{align}
I_{n}^\phi + \sum_{j \in \cN_{\rightarrow n}^\phi} \hspace{-.2cm} I_{jn}^\phi  - \sum_{k \in \cN_{n \rightarrow}^\phi}  \hspace{-.2cm} \, I_{nk}^\phi  = 0
\label{kcl}
\end{align}
where $\cN_{\rightarrow n}^\phi := \{j : (j,n) \in \cE, \phi \in \cP_n \cap \cP_{jn} \}$, and $\cN_{n \rightarrow}^\phi := \{k : (n,k) \in \cE, \phi \in \cP_n \cap \cP_{nk}\}$.  

Let $S_{L_n}^{\phi} := P_{L_n}^{\phi} + j Q_{L_n}^{\phi}$ denote the conglomerate load demanded by residential and commercial facilities
on phase $\phi$ of node $n$, and $S_{G_n}^{\phi} := P_{G_n}^{\phi} + j Q_{G_n}^{\phi}$ the overall power supplied by conventional distributed generation (DG) units, if any. Loads and distributed generators are modeled as constant PQ units. Suppose further that $R_n^\phi$ RESs  (e.g., PV systems, small wind turbines, or a combination of both) are installed at the same phase and node, and let $S_{E_n,r}^{\phi} := P_{E_n,r}^{\phi} + j Q_{E_n,r}^{\phi}$ denote the \emph{actual} power supplied by RES $r$.\footnote{If DG units and RESs are required to operate at unitary power factor, their supplied reactive power is set to zero; that is, $Q_{G_n}^{\phi} = Q_{E_n,r}^{\phi} = 0$. }  
Overall, power balance at phase $\phi \in \cP_n$ of node $n \in \cN$ implies that    
\begin{align}
\label{loadmodel}
V_n^{\phi} (I_n^{\phi})^* =   S_{G_n}^{\phi} + \sum_{r = 1}^{R_n^{\phi}} S_{E_n,r}^{\phi} -  S_{L_n}^{\phi} \, .
\end{align}

With or without uncertainty present, the objective of distribution system reconfiguration schemes is to identify the network topology that is optimal in a well-defined sense, while ensuring load demand satisfaction and adherance to thermal and security constraints~\cite{Khodr09,Baran89,Schmidt95,Cicoria04}. Traditionally, the sought configuration is the radial one~\cite{Khodr09,Baran89,Schmidt95}, although meshed networks were also explored in~\cite{Cicoria04,Dallanese-TPD13}. Unfortunately, variables $\{S_{L_n}^{\phi}\}$  are generally affected by load forecasting errors, whereas the instantaneous power $\{S_{E_n,r}^{\phi}\}$ harvested by RESs will conceivably fluctuate around its forecasted values due to e.g., fast-varying weather conditions~\cite{Bacher2009,Tsikalakis06}. Thus, it is essential to account for possible supply-demand imbalance emerging from uncertain RES generation and load forecasting errors, in order to ensure a reliable power delivery microgrid-wide, and make risk-limiting operational decisions~\cite{Celli02}. 

The goal here is a microgrid configuration that ensures \emph{satisfaction of the load} with arbitrarily high probability, while at the same time adhering to thermal constraints. To this end, notice first that a loss of load occurs whenever the net power $V_n^{\phi} (I_n^{\phi})^* - S_{G_n}^{\phi} - \sum_{r} S_{E_n,r}^{\phi} $ exiting node $n$ and phase $\phi$ is not sufficient to satisfy the load; that is, when $- S_{L_n}^{\phi} < V_n^{\phi} (I_n^{\phi})^* -  S_{G_n}^{\phi} - \sum_{r} S_{E_n,r}^{\phi}$. Let $\cD := \{(\phi, n): S_{E_n,r}^{\phi} \neq 0, \mathrm{~or~}  S_{L_n}^{\phi} \neq 0\}$ denote the set collecting the phase-node pairs where generators and/or loads are located, and define $\bar{\cD} := \{(\phi, n): S_{E_n,r}^{\phi} =  S_{L_n}^{\phi} = 0\}$. Then, with $\rho \in (0,1)$ representing a pre-selected threshold for the LOL probability, and upon defining the vector-valued function 
\begin{align}
\hspace{-.2cm} \bell_{n}^{\phi}(\cV) :=  
\left[ 
\begin{array}{ll}
 \Re\{V_n^{\phi} (I_n^{\phi})^*\} - P_{G_n}^{\phi} - \sum_{r} P_{E_n,r}^{\phi} +  P_{L_n}^{\phi}  \\
 \Im\{V_n^{\phi} (I_n^{\phi})^*\} - Q_{G_n}^{\phi} - \sum_{r} Q_{E_n,r}^{\phi} +  Q_{L_n}^{\phi}  
\end{array}
\right] \nonumber 
\end{align}
\hspace{-.2cm} where $\cV := \{\{x_{mn}\}, \{I_{mn}^\phi\}, \{I_n^\phi, V_n^\phi\}, \{P_{G_n}^{\phi},Q_{G_n}^{\phi}\} \}$ collects the microgrid design variables, the following constraint enforces \emph{every} load to be satisfied with probability at least $1 - \rho$:
\begin{align}
\label{solprob}
\Prob\left\{\bell_{n}^{\phi}(\cV)  \preceq \mathbf{0}, \forall \, (\phi, n) \in \cD  \right\} \geq 1 - \rho \, .
\end{align}

Based on~\eqref{solprob}, the novel risk-constrained microgrid reconfiguration task can be formulated as:
\begin{subequations}
\label{recon_nonconvex_current}
\begin{align}
& \hspace{-1.8cm} \mathrm{(MR1)}   \quad \quad \min_{\cV} \quad  C(\cV)  \label{P1cost}  \\
\textrm{subject~to~~}  & \eqref{line_selection},~\eqref{kcl},~\eqref{loadmodel},~\eqref{solprob} ~\mathrm{and} \nonumber \\
& \hspace{-.35cm} I_n^{\phi} = 0 \, ,\hspace{1.55cm}  \forall \, (\phi, n) \in \bar{\cD} \\
& \hspace{-.75cm}  |I_{mn}^{\phi}| \leq I_{mn}^{\textrm{max}},    \hspace{1.15cm} \forall  \, \phi \in \cP_{mn},  (m,n) \in \cE \label{P1magnitude} \\
& \hspace{-.6cm} x_{mn}  \in \{0,1\} ,    \, \hspace{.9cm}  \forall \, (m,n) \in \cE_R  \hspace{.9cm}   \label{P1binary} \\
& \hspace{-.6cm}  S_n^\textrm{min} \preceq S_n^\phi \preceq S_n^\textrm{max} \, , \hspace{.2cm}  \forall \, \textrm{DG unit}   \label{P1dglimits}
%& \hspace{-.4cm} V_n^\phi \in \cB_n^{\phi} \, , \hspace{1.3cm}  \phi \in \cP_n, \ n \in \cN   
\end{align}
\end{subequations}
where $C(\cV)$ is given cost; $I_{mn}^{\textrm{max}}$ is a cap for $|I_{mn}^{\phi}|$
to protect conductors from overheating; and,~\eqref{P1dglimits} are box constraints for the DG units. When the objective is to minimize the overall active power loss~\cite{Baran89,Schmidt95}, the cost is selected to be
%\begin{align}
$C(\cV) = \sum_{(m,n) \in \cE} \Re\{ \bi_{mn}^{\cT} \bZ_{mn} \bi_{mn}\}$. 
%\label{power_loss}
%\end{align}
Alternatively, the net microgrid operation cost can be minimized by setting $C(\cV) = \sum_{\phi \{a,b,c\}} \hspace{-.2cm} c_1 \Re\{V_1^\phi (I_1^\phi)^{*}\} +  \sum_{n \in \cN \backslash \{1\} ,\phi \in \cP_n} \hspace{-.2cm} c_n^\phi P_{G,n}$, with $c_1$ and $c_n^\phi$ representing the costs of power drawn at the PCC and supplied by the conventional DG at node $n$ and phase $\phi$, respectively. A weighted combination of the two can also be employed, along with (convex) terms (e.g., $\alpha_{mn} |I_{mn}^{\phi}|$) to account for possible line maintenance and security costs~\cite{Celli02}. In the presence of dispatchable loads, a disutility function can be introduced to capture end-user dissatisfaction when operating away from a nominal point. Furthermore, an optimization variable can also be included in the balance equation~\eqref{solprob} to account for the amount of load curtailed.  

To appreciate the value of such problem formulation, notice that (MR1) an be employed to to reconfigure the microgrid after a localized outage, without accurate information on load and RES generation. As for operation planning, (MR1) can be useful to decide whether or not the microgrid can afford operating in an islanded mode without incurring LOL.      

Unfortunately, solving (MR1) is computationally prohibitive for three reasons: 

\noindent \emph{r1)} due to the binary variables $\{x_{mn}\}$, solving (MR1) is \emph{NP-hard}; finding the globally optimal set of binary variables requires solving $2^{|\cE_{R}|}$ subproblems;

\noindent \emph{r2)} the bilinear terms $x_{mn} \bv_m$  and $V_n^{\phi} (I_n^{\phi})^*$ in~\eqref{line_selection} and~\eqref{solprob}, respectively, render (MR1) \emph{nonconvex}; even for fixed values of $\{x_{mn}\}$, nonconvexity implies that (MR1) is difficult to solve optimally and efficiently; and,  

\noindent \emph{r3)} the probabilistic constraint~\eqref{solprob} is generally in a computationally intractable form. To obtain a tractable surrogate constraint, it is first necessary to find the probability distribution function (pdf) of the random variables $\{ \sum_{r =1}^{R_n^\phi} S_{E_n,r}^{\phi} -  S_{L_n}^{\phi} \}$. This is however, a major challenge on its own. In fact, while for single wind farm or PV system ($R_n^\phi = 1$)
this is possible~\cite{Bacher2009,Liu2010}, the pdf of the power supplied by multiple wind farms and PV systems (along with the load) is hard to obtain. And, even if a pdf becomes available, it may not lead to a convex re-formulation of~\eqref{solprob}.  

One approach to coping with \emph{r1)}--\emph{r3)} is proposed in the ensuing section, along with a computationally tractable re-formulation of (MR1).  

%%%%%%%%%%%%%%%%%%%%%%%%%%%%%%%%%%%%%%%%%
\section{Computationally tractable formulation} 
\label{sec:dsr_sparsity}
%%%%%%%%%%%%%%%%%%%%%%%%%%%%%%%%%%%%%%%%%

Collect first the real and imaginary parts of $I_{n}^\phi$ in the $2|\cP_{n}| \times 1$ vector $\biota_{n}^\phi := [\Re^\cT\{I_{n}^\phi\}, \Im^\cT\{I_{n}^\phi\}]^{\cT} \in \mathbb{R}^{2}$; and likewise define the vector $\bxi_{mn} := [\Re^\cT\{\bi_{mn}\}, \Im^\cT\{\bi_{mn}\}]^{\cT} \in \mathbb{R}^{2|\cP_{mn}|}$. % and stack in the real-valued $2 \times 1$ vectors $\bsigma_{G_n}^\phi$, $\bsigma_{L_n}^\phi$, and $\bsigma_{E_n,r}^\phi$ the real and imaginary parts of $S_{G_n}^\phi$, $S_{L_n}^\phi$, and $S_{E_n,r}^\phi$, respectively. 

To bypass \emph{r1)}, the approach in~\cite{Dallanese-TPD13} is broadened here to account for load and RES generation uncertainty. To this end, notice first that the entries of $\bxi_{mn}$  are \emph{all} zero if line $(m,n) \in \cE_R$ is \emph{not} used to deliver power to the loads; that is, $I_{mn}^{\phi} = 0$ for all phases $\phi \in \cP_{mn}$. Clearly, $\bxi_{mn} \neq \mathbf{0}$ otherwise. Adopting the compressive sampling terminology~\cite{YuLi06}, vector $\bxi_R := [\{\bxi_{mn} | (m,n) \in \cE_R\}]^{\cT}$ is \emph{group sparse}, meaning that ``group(s) of elements'' ($\bxi_{mn}$ in this case) are either all zero, or not. One major implication of this group sparsity attribute of $\bxi_R$, is that one can discard the binary variables $\{x_{mn}\}_{(m,n) \in \cE_R}$, and effect line selection by augmenting the cost~\eqref{P1cost} with the following convex group-Lasso-type regularization term~\cite{YuLi06,Wiesel11}
\begin{align}
g(\bxi_R) :=  \lambda \hspace{-.3cm} \sum_{(m,n) \in \cE_{R}} \, \| \bxi_{mn} \|_2  \label{Glasso_currents}  
\end{align}
where $\lambda > 0$ is a tuning parameter. Specifically, the role of $\lambda$ is to control the number of vectors $\{\bxi_{mn}\}_{(m,n) \in \cE_R}$ (and, hence currents $\{\bi_{mn}\}_{(m,n) \in \cE_R}$) that are set to zero. This means that by adjusting $\lambda$ one can obtain either meshed topologies (low values of $\lambda$), weakly-meshed, or even radial systems (high values of $\lambda$); see e.g.~\cite{Wiesel11,Dallanese-TPD13}. A generalization of~\eqref{Glasso_currents} is represented by the weighted version $g_w(\bxi_R) :=  \sum_{(m,n) \in \cE_{R}} \, \lambda_{mn} \| \bxi_{mn} \|_2 $, where $\{\lambda_{mn}\}$ substantiate possible operator preferences to use (low value of $\lambda_{mn}$) or not (high value of $\lambda_{mn}$) specific lines. For instance, it may be preferable to open or close switches that are commanded remotely, rather than requiring hand operations in situ.  

To avoid the bilinear terms $\{V_n^\phi (I_n^\phi)\}$, and considerably lower the complexity incurred by the resultant optimization scheme,  
consider adopting the approximate current-power relation employed by~\cite{Bolognani13}. 
Specifically, with $V_N = M_N e^{\varphi_N^\phi}$ denoting the nominal line-to-ground voltage on phase $\phi$, the injected current can be approximated as $I_n^\phi \approx (1/M_N)e^{j \varphi_N^\phi} (S_{G_n}^{\phi} + \sum_{r} S_{E_n,r}^{\phi} -  S_{L_n}^{\phi})^*$. Although this approach provides a surrogate \emph{linear} (as opposed to bilinear) load balance equation, the approximation error introduced must be carefully accounted for in~\eqref{solprob}. To this end, the \emph{actual} current injected at phase $\phi$ of node $n$ is modeled here as   
\begin{align}
\label{approx_current}
\biota_n^\phi  :=  \bPhi_n^\phi  \left(
\left[ \hspace{-.2cm} 
\begin{array}{rr}
P_{G_n}^{\phi} \\
Q_{G_n}^{\phi} 
\end{array} \hspace{-.2cm}
\right] +
\sum_{r} 
\left[ \hspace{-.2cm} 
\begin{array}{rr}
P_{E_n}^{\phi} \\
Q_{E_n}^{\phi} 
\end{array} \hspace{-.2cm}
\right]
- \left[ \hspace{-.2cm} 
\begin{array}{rr}
P_{L_n}^{\phi} \\
Q_{L_n}^{\phi} 
\end{array} \hspace{-.2cm}
\right] \right) + \bepsilon_{\iota_n}^\phi 
\end{align}
where $\bepsilon_{\iota_n}^\phi$ captures approximation errors, and $\bPhi_n^\phi$ is a $2 \times 2$ matrix with columns
%\begin{align}
%\bPhi_n^\phi := \frac{1}{M_N} 
%\left[ 
%\begin{array}{rr}
%\Re\{e^{\varphi_N^\phi}\} & \Im\{e^{\varphi_N^\phi}\} \\
%\Im\{e^{\varphi_N^\phi}\} & - \Re\{e^{\varphi_N^\phi}\} 
%\end{array}
%\right] . \nonumber 
%\end{align} 
$(1/M_N)[ \Re\{e^{\varphi_N^\phi}\} , \Im\{e^{\varphi_N^\phi}\}]^\cT$ and $(1/M_N) [\Im\{e^{\varphi_N^\phi}\} , - \Re\{e^{\varphi_N^\phi}\}]^\cT$. Define now the vector function
\begin{subequations}
\begin{align}
\tilde{\bell}_{n}^{\phi} :=  
\left[ \hspace{-.2cm}
\begin{array}{ll}
(\bvarphi_{n}^{\phi})^\cT \biota_n^\phi  & - P_{G_n}^\phi -  \sum_{r} P_{E_n,r}^\phi + P_{L_n}^\phi  - (\bvarphi_{n}^{\phi})^\cT  \bepsilon_{\iota_n} \\
(\bar{\bvarphi}_{n}^{\phi})^\cT \biota_n^\phi & -  Q_{L_n}^\phi -  \sum_{r} Q_{E_n,r}^\phi + Q_{L_n}^\phi  - (\bar{\bvarphi}_{n}^{\phi})^\cT  \bepsilon_{\iota_n}  
\end{array}\hspace{-.2cm} \right] \nonumber
\end{align}  
\end{subequations}
where $\bvarphi_{n}^{\phi} := [ M_N \Re\{e^{\varphi_N^\phi}\} , M_N \Im\{e^{\varphi_N^\phi}\}]^\cT$ and $\bar{\bvarphi}_{n}^{\phi} := [M_N \Im\{e^{\varphi_N^\phi}\} , - M_N \Re\{e^{\varphi_N^\phi}\}]^\cT$. Then, constraint~\eqref{solprob} can be equivalently re-written as  
\begin{align}
\hspace{-.2cm} \Prob \Big\{\tilde{\bell}_{n}^{\phi}(\bxi,\{S_n^\phi\}) \preceq \mathbf{0} , \quad  \forall (\phi, n) \in \cD \Big\} \geq 1 - \rho
 \label{lolprob_2}
\end{align}  
where the probability is evaluated over the pdf of random variables $\{\sum P_{E_n,r}^\phi - P_{L_n}^\phi + (\bvarphi_{n}^{\phi})^\cT  \bepsilon_{\iota_n}\}, \{\sum Q_{E_n,r}^\phi - Q_{L_n}^\phi + (\bar{\bvarphi}_{n}^{\phi})^\cT  \bepsilon_{\iota_n}\}$. The empirical pdf of $\bepsilon_{\iota_n}^\phi$ can be obtained either using historical data, or, from the voltage distribution~\cite{Bolognani13}. 

Consider re-expressing the current $\biota_n^\phi$ as $\biota_n^\phi = \bmA_n^\phi \bxi$, where $\bxi$ stacks all the line current vectors $\{\bxi_{mn}\}$, and $\bmA_n^\phi$ is obtained in the obvious way from Kirchhoff's current law~\eqref{kcl}. Then, based on~\eqref{Glasso_currents} and~\eqref{lolprob_2}, the microgrid reconfiguration problem can be reformulated as: 
\begin{subequations}
\label{recon_infinite_dimensional}
\begin{align}
& \hspace{-.7cm} \mathrm{(MR2)}   \quad  \min_{\bxi, \{S_n^\textrm{min} \preceq S_n^\phi \preceq S_n^\textrm{max}\}} \quad  C(\bxi, \{S_{G_n}^\phi\}) + g(\bxi_R) \label{P2cost}  \\
\textrm{s.t.~~}  &  \hspace{-.0cm}  \bxi_{mn}^{\cT} \bbmM_{mn}^{\phi} \bxi_{mn} \leq (I_{mn}^{\textrm{max}})^2 ,  \,\,  \, \phi \in \cP_{mn},  (m,n) \in \cE \label{P2magnitude} \\
& \hspace{.1cm} \Prob \Big\{\tilde{\bell}_{n}^{\phi}(\bxi,\{S_n^\phi\}) \preceq \mathbf{0} ,   \quad  \forall (\phi,n) \in \cD \Big\} \geq 1 - \rho \label{P2solp}  \\
& \hspace{.2cm}  \bmA_n^\phi \bxi  = \mathbf{0}_{2 \times 1} \, , \forall \, (\phi, n) \in \bar{\cD} 
\end{align}
\end{subequations}
where $\bbmM_{mn}^{\phi} := \bI_{2} \otimes \be_{mn}^\phi (\be_{mn}^\phi)^\cT$, with $\{\be_{mn}^{\phi}\}_{\phi \in \cP_{mn}}$ representing the canonical basis of $\mathbb{R}^{|\cP_{mn}|}$. 
%One way to transform~\eqref{P2solp} in a tractable form is to approximate $\sum \bsigma_{E_n,r}^\phi - \bsigma_{L_n}^\phi + (\bPhi_n^\phi)^{-1}  \bepsilon_{\iota_n}^\phi$ as Gaussian distributed (for sufficiently large values of $R_n^\phi$), by appealing to the central limit theorem. Unfortunately, central limit theorem relies on vanishing correlation among the random summands,  and this does not hold in practice for $\{\bsigma_{E_n,r}^\phi\}$ (especially for geographically close wind turbines and PV systems). 
To address \emph{r3)}, a computationally efficient scheme is presented next, based on the so-called scenario-based convex approximation~\cite{Calafiore}.

\vspace{-.3cm}

%%%%%%%%%%%%%%%%%%%%%%%%
\subsection{Scenario-based approximation } 
\label{sec:scenario}
%%%%%%%%%%%%%%%%%%%%%%%%

To briefly illustrate the general scenario-based approximation method~\cite{Calafiore}, consider the prototype convex problem:
\begin{subequations}
\label{scenario_approximation}
\begin{align}
& \hspace{-1cm} \mathrm{(P)}   \quad  \min_{\bx \in \cX} \quad  c(\bx)  \label{Pcost}  \\
\textrm{s.t.~~}  & \Prob\{ \bff(\bx,\bsigma) \preceq \mathbf{0}  \}  \geq 1 - \rho \label{Psolp}  
\end{align}
\end{subequations}
where $\cX \subseteq \mathbb{R}^m$ is a nonempty convex set; $c: \cX \rightarrow \mathbb{R}$ is convex; $\bsigma$ is a random vector, whose pdf has support $\cU \subseteq \mathbb{R}^n$; and, $\bff: \cX \times \cU \rightarrow \mathbb{R}^d$ is a vector-valued convex function.  Then, the scenario-based approximation method amounts to: \emph{i)} generating $K$ independent samples $\bsigma(1), \ldots, \bsigma(K)$; and, \emph{ii)} approximating (P) as the following convex program  
\begin{subequations}
\label{scenario_approximation2}
\begin{align}
& \hspace{-1cm} \mathrm{(PA)}   \quad  \min_{\bx \in \cX} \quad  c(\bx)  \label{PAcost}  \\
\textrm{s.t.~~}  & \bff(\bx,\bsigma(k)) \preceq \mathbf{0} \, , \forall j = 1,\ldots, K .  \label{PAsolp}  
\end{align}
\end{subequations}  
To better appreciate the merits of this approach, notice first that since $c(\cdot)$ and $\bff(\cdot)$ are convex, (PA) is a \emph{convex} program. Further, to derive the approximate (PA), no specific requirements on the distribution of $\bsigma$ are imposed. However, a pertinent question is whether the solution of (PA) is feasible also for the original problem (P), given that the constraints $\bff(\bx,\bsigma(k)) \preceq \mathbf{0}$ are randomly selected, and the resulting optimal solution $\bx^{(PA)}$
is a random variable that depends on the extracted samples $\bsigma(1), \ldots, \bsigma(K)$. Let $\beta$ denote a cap for the probability of $\bx^{(PA)}$ being not feasible for (P) (also referred to as the ``risk of failure''~\cite{Calafiore}). Then, given $\rho$ and $\beta$, it can be shown that if the number of samples $K$ is chosen such that (see~\cite[Corollary~1]{Calafiore})
\begin{align}
K \geq \tilde{K} := \ceil[\Big]{2 \rho^{-1} \ln \beta^{-1}  + 2 m + 2 m \rho^{-1}  \ln 2 \rho^{-1} }
\label{sample_size}
\end{align}
then the optimal solution to (PA) is feasible for (P) with probability at least $1-\beta$.

To apply the scenario-based approximation method to (MR2), generate $K$ independent samples $P_{n}^\phi(k) := \sum_{r} P_{E_n,r}^\phi(k) - P_{L_n}^\phi(k) + (\bvarphi_{n}^{\phi})^\cT  \bepsilon_{\iota_n}(k)$, and $Q_{n}^\phi(k) := \sum_{r} Q_{E_n,r}^\phi(k) - Q_{L_n}^\phi(k) + (\bvarphi_{n}^{\phi})^\cT  \bepsilon_{\iota_n}(k)$, $k = 1,\ldots,K$, and replace the chance-constraint~\eqref{P2solp} with the linear constraints 
\begin{subequations}
\label{approx_solp}
\begin{align}
& (\bvarphi_{n}^{\phi})^\cT \bmA_n^\phi \bxi  - P_{G_n}^\phi  \leq P_{n}^\phi(k) \, , \, k = 1, \ldots, K     \\
& (\bar{\bvarphi}_{n}^{\phi})^\cT \bmA_n^\phi \bxi  - Q_{G_n}^\phi \leq Q_{n}^\phi(k) \, , \, k = 1, \ldots, K \, .\hspace{-.1cm}   
\end{align}
\end{subequations}
One possible limitation of this approach is that the minimum number of samples $\tilde{K}$ increases rapidly as $\rho$ decreases. Further, $\tilde{K}$ is very large for microgrids of medium- large-size (in fact, the total number of constraints would amount to $\tilde{K} |\cD|$). Luckily, a closer look to~\eqref{approx_solp} reveals that~\eqref{P2solp} can be replaced by the following $2 |\cD|$ constraints 
\begin{subequations}
\label{approx_solp2}
\begin{align}
(\bvarphi_{n}^{\phi})^\cT \bmA_n^\phi \bxi  &  \leq P_{G_n}^\phi + \min_{k = 1, \ldots, K} \{ P_{n}^\phi(k) \}  \, , \forall \, (\phi, n) \in \cD \\
(\bar{\bvarphi}_{n}^{\phi})^\cT \bmA_n^\phi \bxi  &  \leq Q_{G_n}^\phi  + \min_{k = 1, \ldots, K} \{Q_{n}^\phi(k) \}  \, , \forall \, (\phi, n) \in \cD . \hspace{-.2cm} 
\end{align}
\end{subequations}
Thus, replacing~\eqref{P2solp} with~\eqref{approx_solp2}, the following surrogate problem is readily obtained
\begin{subequations}
\label{recon_finite_dimensional}
\begin{align}
& \hspace{-.7cm} \mathrm{(MR3)}   \quad \quad \min_{\bxi, \{\bsigma_{G_n}^\phi\}} \quad  C(\bxi, \{\bsigma_{G_n}^\phi\}) + g(\bxi_R) \label{P3cost}  \\
\textrm{s.t.~~}  & \eqref{approx_solp2},\textrm{~and} \nonumber   \\
& \hspace{-.0cm}  \bxi_{mn}^{\cT} \bbmM_{mn}^{\phi} \bxi_{mn} \leq (I_{mn}^{\textrm{max}})^2 ,  \,\,  \, \phi \in \cP_{mn},  (m,n) \in \cE \label{P3magnitude} \\
& \hspace{1.4cm}  \bmA_n^\phi \bxi  = \mathbf{0}_{2 \times 1} \, , \forall \, (\phi, n) \in \bar{\cD} \, . \label{P3zeroCurrent}
\end{align}
\end{subequations}
If $C(\bxi, \{\bsigma_{G_n}^\phi\})$ is chosen convex as in~\cite{Khodr09,Celli02,Baran89,Schmidt95,Cicoria04}, then~\eqref{recon_finite_dimensional} is a \emph{convex} program that can be solved efficiently via standard interior-point methods. Further, tailoring~\eqref{sample_size} to (MR3), the minimum sample size $\tilde{K}$ is established in the following proposition. 

\vspace{.1cm}

\begin{proposition}
\label{prop:sample_size}
Given the LOL probability threshold $\rho$, and the lower bound 
\begin{align}
& K \geq \tilde{K}^{\textrm{MR3}} := \ceil[\Big]{2 \rho^{-1} \ln \beta^{-1}  + 4 (N_G + \sum_{(m,n) \in \cE} |\cP_{mn}|) \nonumber \\ 
& \hspace{2cm} + 4 (N_G + \sum_{(m,n) \in \cE} |\cP_{mn}|) \rho^{-1}  \ln 2 \rho^{-1} }
\label{sample_size_MR3}
\end{align}
where $N_G$ stands for the total number of conventional DG units, then the solution $(\bxi^{opt}, \{S_{G_n}^{\phi , opt}\})$ to (MR3) is feasible for (MR2) with probability no less than $1 - \beta$. \hfill $\Box$
\end{proposition}

\vspace{.1cm}

Once (MR3) is solved, the optimal topology of the microgrid is obtained by discarding the distribution lines with an associated zero current; that is, 
%\begin{align}
$\cE^{opt} := \cE \backslash \{(m,n) \in \cE_R : \bxi_{mn}^{opt} = \mathbf{0} \}$.
%\label{optimal_topology}
%\end{align}
Given the optimal configuration $(\cN,\cE^{opt})$, voltages and currents can be computed via OPF in real-time, once RES generation and load are revealed; see e.g.,~\cite{Dallanese-TSG13} and references therein. Finally, notice that voltage constraints can be included in (MR3) as described in~\cite{Dallanese-TPD13}.

\vspace{.1cm}

\noindent \emph{Remark 1}. It is worth mentioning that one can alternatively solve (MR1)  using the scenario-based approximation method in conjunction with off-the-shelf solvers for mixed-integer nonlinear programs (MINP). However, this may not be as convenient for three reasons: \emph{i)} for a given sample set $\{\sum_r \bsigma_{E_n,r}^{\phi}(k) - \bsigma_{L_n}^{\phi}(k)\}$, (MR1) is a \emph{nonconvex} program; as a  consequence,~\eqref{sample_size} may not hold, since it is grounded on a convexity assumption~\cite{Calafiore}; \emph{ii)} the computational burden of MINP solvers is typically much higher then that of interior-point methods for convex programs; and, \emph{iii)} solving an MINP in a distributed fashion is not immediate. In contrast, a distributed solver for (MR3) is feasible as shown in Section~\ref{sec:Distributed}.  

\vspace{.1cm}

\emph{Remark 2}. Since values of the diagonal entries of matrix $\bY_{mn}$ are typically on the order of $10-100$ micro Siemens per mile (see e.g.,~\cite{testfeeder}), the effect of line shunt admittances was neglected in the Kirchhoff's current law~\eqref{kcl} (see also~\cite{Khodr09,Celli02,Baran89,Schmidt95,Bolognani13}). Simulation results on IEEE test feeders~\cite{testfeeder} showed that the approximation error is negligible, meaning that the resultant optimal topology does not change upon considering $\bY_{mn}$. However, to account for possible perceptible effects of these shunt admittances in other real-world distribution systems, the shunt elements $(1/2) [\bY_{mn}]_{\phi,\phi}$, $\phi \in \cP_{mn}$, of a line $(m,n)$ can be taken to be constant-admittance loads at nodes $m$ and $n$~\cite{LoadModel}. Thus, upon computing the overall constant-impedance load per node, one can readily obtain an approximate value of the current absorbed by this constant-impedance loads as shown in~\cite{Bolognani13,LoadModel}, and add this approximate current in~\eqref{kcl}. 

%%%%%%%%%%%%%%%%%%%%%%%
\subsection{Multi-period optimization} 
\label{sec:multiperiod}
%%%%%%%%%%%%%%%%%%%%%%%

The system reconfiguration problem (MR1) can be extended to accommodate energy storage systems. To this end, consider optimizing the operation of a distribution system over a (rolling) horizon $t = 1,\ldots,T$~\cite{Palma-Behnke13}, with granularity that depends on whether the ambient conditions are fast-, slow-changing or invariant. Throughout this subsection, let the superscript $(\cdot)^{t}$ index the time slots, and let $B_{n}^{\phi,t} $ represent the state of charge of an energy storage unit located on node $n$ and phase $\phi$ at slot $t$. Then, a multi-period system reconfiguration problem can be formulated a follows: 
\begin{subequations}
\label{recon_nonconvex_current_multip}
\begin{align}
& \hspace{-1.8cm}  \quad \quad \min_{\{\cV^t\}_{t = 1}^T} \quad  \sum_{t = 1}^T C^t(\cV^t)  \label{P1cost_multip}  \\
\textrm{subject~to~~}  & \eqref{line_selection},~\eqref{kcl},~\mathrm{and} \nonumber \\
%& \hspace{-1.45cm}  V_n^{\phi,t} (I_n^{\phi,t})^* =   S_{G_n}^{\phi,t} + \sum_{r} S_{E_n,r}^{\phi,t} -  S_{L_n}^{\phi,t} + P_{B_n}^{\phi,t}, \nonumber \\
%& \hspace{1.8cm}  \forall \, (\phi, n) \in \cD, \forall \, t \\
& \hspace{-.75cm} I_n^{\phi,t} = 0 \, ,\hspace{1.15cm}  \forall \, (\phi, n) \in \bar{\cD}, \forall \, t  \\
& \hspace{-.95cm}  |I_{mn}^{\phi,t}| \leq I_{mn}^{\textrm{max}},    \hspace{.75cm} \forall  \, \phi \in \cP_{mn}, \forall \, (m,n) \in \cE , \forall \, t \label{P1magnitude_multip} \\
& \hspace{-.8cm} x_{mn}^t  \in \{0,1\} ,    \, \hspace{.5cm}  \forall \, (m,n) \in \cE_R , \forall \, t  \hspace{.9cm}   \label{P1binary_multip} \\
& \hspace{-.8cm}  S_n^\textrm{min} \preceq S_n^{\phi,t} \preceq S_n^\textrm{max} \, , \hspace{.2cm}  \forall \, \textrm{DG unit} , \forall \, t   \label{P1dglimits_multip}\\
& \hspace{-.95cm} 1 - \rho \leq \Prob\left\{\bell_{n}^{\phi,t}(\cV^t)  \preceq \mathbf{0}, \forall \, (\phi, n) \in \cD, \forall t  \right\}  \label{LOL_multip3}  \\
& \hspace{-.95cm} B_{n}^{\phi,t+1} = B_{n}^{\phi,t} + P_{B_n}^{\phi,t} \hspace{1.1cm}  \forall \, t = 1,\ldots,T-1  \label{P1battery_multip} \\
& \hspace{-.95cm} B_{n}^{\phi,\mathrm{min}} \leq B_{n}^{\phi,t+1}  \leq B_{n}^{\phi,\mathrm{max}} \hspace{.3cm}  \forall \, t = 1,\ldots,T-1 \label{P1battery_multip2} \\
& \hspace{-1.8cm}  - \eta_n^{\phi,\mathrm{dis}} B_{n}^{\phi,t}  \leq P_{B_n}^{\phi,t}  \leq \eta_h^{\mathrm{ch}} (B_{n}^{\mathrm{max}} - B_{n}^{\phi,t} )  \hspace{.2cm}  \forall \, t \label{P1battery_multip3} 
\end{align}
\end{subequations}
where constraints~\eqref{line_selection} and~\eqref{kcl} are enforced per slot $t = 1,\ldots T$; 
$\cV^t := \{\{x_{mn}^t\}, \{I_{mn}^{\phi,t}\}, \{I_n^{\phi,t}, V_n^{\phi,t}\}, \{P_{G_n}^{\phi,t},Q_{G_n}^{\phi,t}\}, \{P_{B_n}^{\phi,t}\} \}$ collects the optimization variables pertaining to time slot $t$;~\eqref{P1battery_multip} is the dynamical equation of the energy storage system; $\eta_n^{\phi,\mathrm{ch}}$ and $\eta_h^{\phi,\mathrm{dis}}$ are charging and discharging efficiencies; and, $\bell_{n}^{\phi,t}(\cV^t)$ is re-defined as: 

{\small
\begin{align}
 \bell_{n}^{\phi}(\cV) :=  
\left[ 
\begin{array}{ll}
 \Re\{V_n^{\phi,t} (I_n^{\phi,t})^*\} - P_{G_n}^{\phi,t} - \sum_{r} P_{E_n,r}^{\phi,t} +  P_{L_n}^{\phi,t} +  P_{B_n}^{\phi,t}\\
 \Im\{V_n^{\phi,t} (I_n^{\phi^t})^*\} - Q_{G_n}^{\phi,t} - \sum_{r} Q_{E_n,r}^{\phi,t} +  Q_{L_n}^{\phi,t}  
\end{array}
\right] \nonumber 
\end{align}}
\normalsize

Towards obtaining a convex relaxation of problem~\eqref{recon_nonconvex_current_multip}, group-Lasso-type regularization terms can be used to bypass binary selection variables. Specifically, the regularization function 
$g(\{\bxi_R^t\}) :=  \sum_{t=1}^T \lambda^t \sum_{(m,n) \in \cE_{R}} \, \| \bxi_{mn}^t \|_2 $ effects line selection at each slot $t$. 
The scenario-based convex approximation outlined in Section~\ref{sec:scenario} can be utilized to obtain an approximate yet tractable reformulation of constraint~\eqref{LOL_multip3}. The battery (dis)charging can be fine-tuned during the real-time system operation using tools such as OPF.
 
%%%%%%%%%%%%%%%%%%%%%%%%
%\subsection{LOLP vis-\`a-vis worst-case} 
%\label{sec:robust}
%%%%%%%%%%%%%%%%%%%%%%%%

%%%%%%%%%%%%%%%%%%%%%%%%%%%%%%%%%%%%%%%%%%%%%%
\section{Distributed algorithm}
\label{sec:Distributed}
%%%%%%%%%%%%%%%%%%%%%%%%%%%%%%%%%%%%%%%%%%%%%%

A distributed reconfiguration algorithm is desirable when the microgrid includes single- or multi-facility clusters  that are managed independently from the rest of the network in order to pursue individual economic interests. Thus, each cluster autonomously selects the topology of its own subnetwork, and controls the power supplied by conventional DG units and RESs (e.g., for reactive compensation~\cite{Bolognani13}). 
Consider  partitioning the microgrid into $L$ areas $\{\cA^{(\ell)} \subset \cN\}_{\ell = 1}^{L}$, and let $\cN^{(\ell)} := \{j | \exists (m,n) \in \cE: m \in \cA^{(\ell)}, n \in \cA^{(j)}\}$ denote the set of neighboring areas for the $\ell$-th one. Further, let function $\cI(\ell,j): \cA^{(\ell)} \times \cA^{(j)} \rightarrow \cE$ identify the line(s) connecting areas $\ell$ and $j$. It clearly holds that $\cI(\ell,j) = \emptyset$ if $j \notin \cN^{(\ell)}$. Autonomous areas are managed by local area controllers (LACs)~\cite{Nordman12}, whereas the microgrid manager (MGM)  controls the lines interconnecting areas, and the remaining portion of the microgrid.  
Let $\bxi^{(\ell)}$ collect the real and imaginary parts of line currents of lines within area $\ell$, plus the lines connecting area $\ell$ to its neighbors; that is, $\cE^{(\ell)}:=\{(m,n) \in \cE| m,n \in \cA^{(\ell)}\} \cup \{(m,n) \in \cE | m \in \cA^{(\ell)}, n \in \cA^{(j)}, j \in \cN^{(\ell)}\}$. Define $\cE^{(\ell)}_R := \{(m,n) \in \cE_R: m,n \in \cA^{(\ell)} \}$, and assume for simplicity that all lines $\cI := \{\cI(\ell,j), \ell, j = 1, \ldots L\}$ are equipped with sectionalizing switches. 
Next, let  $\bchi_{\cI(\ell,j)}$ represent a \emph{copy} of the vector collecting the currents flowing on the line $\cI(\ell,j)$ connecting areas $\ell$ and $j$.
Consider now decomposing the cost function~\eqref{P3cost} as
%\begin{align}
%C(\bxi, \{\bsigma_{G_n}^\phi\}) & = \sum_{\ell = 1}^L \Big[ C^{(\ell)}(\bxi^{(\ell)}, \{\bs_{G}^{(\ell)}\}) + \hspace{-.3cm} \sum_{(m,n)\in\cE_R^{(\ell)}} \hspace{-.3cm} \lambda^{(\ell)} \|\bxi\|_2 \Big] \nonumber \\
%& \hspace{-1.5cm} + \sum_{\ell = 1}^L \sum_{j > \ell}^L \Big[ C_{\cI(\ell,j)}(\bchi_{\cI(\ell,j)}) + \lambda_{\cI(\ell,j)} \|\bchi_{\cI(\ell,j)}\|_2 \Big]
%\label{cost_decomposition}
%\end{align}
$C(\bxi, \{\bsigma_{G_n}^\phi\}) = \sum_{\ell = 1}^L \Big[ C^{(\ell)}(\bxi^{(\ell)}, \{\bs_{G}^{(\ell)}\}) + \sum_{(m,n)\in\cE_R^{(\ell)}} \lambda^{(\ell)} \|\bxi\|_2 \Big] + \sum_{\ell = 1}^L \sum_{j > \ell}^L \Big[ C_{\cI(\ell,j)}(\bchi_{\cI(\ell,j)}) + \lambda_{\cI(\ell,j)} \|\bchi_{\cI(\ell,j)}\|_2 \Big]$, 
where $C^{(\ell)}(\cdot)$ stands for the cost associated with area $\ell$; $\bs_{G}^{(\ell)}$ collects the powers injected by DG units located within $\cA^{(\ell)}$; and, $C_{\cI(\ell,j)}(\cdot)$ is the cost associated with line $\cI(\ell,j)$. Thus, 
(MR3) can be equivalently reformulated as 
\begin{subequations}
\label{recon_finite_dimensional2}
\begin{align}
& \hspace{-.4cm}  \,\, \min_{\bxi, \{\bsigma_{G_n}^\phi\}} \quad   \sum_{\ell = 1}^L \Big[ C^{(\ell)}(\bxi^{(\ell)}, \{\bs_{G}^{(\ell)}\}) + \hspace{-.3cm} \sum_{(m,n)\in\cE_R^{(\ell)}} \hspace{-.3cm} \lambda^{(\ell)} \|\bxi_{mn}^{(\ell)}\|_2 \Big] \nonumber \\
& \hspace{.2cm} + \sum_{\ell = 1}^L \sum_{j > \ell}^L \Big[ C_{\cI(\ell,j)}(\bchi_{\cI(\ell,j)}) + \lambda_{\cI(\ell,j)} \|\bchi_{\cI(\ell,j)}\|_2 \Big] \hspace{-.2cm}  \label{P3bcost}  \\
& \textrm{s.t.~~} \quad  \{\bxi^{(\ell)}, \bs_{G}^{(\ell)}\} \in \cF^{(\ell)} \, ,\quad \ell = 1. \ldots L  \label{Pb3area} \\
& \hspace{.6cm}  \bchi_{\cI(\ell,j)} \bbmM_{\cI(\ell,j)}^{\phi} \bchi_{\cI(\ell,j)} \leq (I_{\cI(\ell,j)}^{\textrm{max}})^2 , \forall \, \cI(\ell,j) \in \cI\label{P3bmagnitude} \\
& \hspace{.6cm} \bxi_{\cI(\ell,j)}^{(\ell)} = \bchi_{\cI(\ell,j)} \, ,  \forall \, \cI(\ell,j) \in \cI\label{P3bconsensus} \\
& \hspace{.6cm} \bxi_{\cI(\ell,j)}^{(j)} = \bchi_{\cI(\ell,j)} \, ,  \forall \, \cI(\ell,j) \in \cI\label{P3bconsensus2}
\end{align}
\end{subequations}
where constraints~\eqref{P3bconsensus}--\eqref{P3bconsensus2} enforce MGM, LAC $\ell$, and LAC $j$ to \emph{consent} on the value of the currents on line $\cI(\ell,j)$;  $\cF^{(\ell)}$ denotes the set of variables $(\bxi^{(\ell)}, \{\bs_{G}^{(\ell)}\} )$ satisfying constraints~\eqref{approx_solp2},~\eqref{P3magnitude}, and~\eqref{P3zeroCurrent} in each phase of nodes $n \in \cA^{(\ell)}$; and, $\bbmM_{\cI(\ell,j)}^{\phi}$ is defined in the obvious way [cf.~\eqref{P2magnitude}].

Constraints~\eqref{P3bconsensus}--\eqref{P3bconsensus2} render problems (MR3) and~\eqref{recon_finite_dimensional2}
equivalent; however, the same constraints couple both optimization problems across areas. To enable a decentralized solution, consider introducing $|\cI|$ auxiliary variables $\{\bz_{\cI(\ell,j)}\}$, and replace~\eqref{P3bconsensus}--\eqref{P3bconsensus2} with the following equivalent set of constraints per line $\cI(\ell,j)$:  
\begin{subequations}
\label{consensus_constraints}
\begin{align}
&  \bxi_{\cI(\ell,j)}^{(\ell)} = \bz_{\cI(\ell,j)} \, ,  \bxi_{\cI(\ell,j)}^{(j)} = \bz_{\cI(\ell,j)} \, , \label{consensus_constraints2} \\
&  \mathrm{and~}  \bz_{\cI(\ell,j)} =  \bchi_{\cI(\ell,j)}  \, . 
\label{consensus_constraints3}
\end{align}
\end{subequations}
The idea here is to solve the resultant optimization problem specified by~\eqref{P3bcost}--\eqref{P3bmagnitude} and~\eqref{consensus_constraints} in a decentralized fashion by resorting to the ADMM~\cite[Sec.~3.4]{BeT89}. In principle,~\eqref{P3bcost}--\eqref{P3bmagnitude},~\eqref{consensus_constraints} can be solved via dual (sub-)gradient ascent iterations~\cite{Nedic09}. However, primal averaging is necessary when the dual function is non-differentiable and the step size is fixed, thus resulting in a typically slower convergence than ADMM~\cite{Dallanese-TSG13}. 

Toward this end,  let $\{\bgamma^{(\ell)}_{\cI(\ell,j)}\}$ and $\{\bmu_{\cI(\ell,j)}\}$ be the vector-valued multipliers associated 
with constraints~\eqref{consensus_constraints2} and~\eqref{consensus_constraints3}, respectively,  and 
consider the following partial quadratically-augmented Lagrangian of~\eqref{P3bcost}--\eqref{P3bmagnitude} and~\eqref{consensus_constraints}: 
\begin{align}
& \hspace{-.0cm} \cL(\{\bxi^{(\ell)}, \bs_{G}^{(\ell)}\}, \{\bchi_{\cI(\ell,j)}\}, \{\bz_{\cI(\ell,j)}\}, \bd) \nonumber \\
& \hspace{-.0cm}  := \sum_{\ell = 1}^L \Big[ C^{(\ell)}(\bxi^{(\ell)}, \{\bs_{G}^{(\ell)}\}) + \hspace{-.3cm} \sum_{(m,n)\in\cE_R^{(\ell)}} \hspace{-.3cm} \lambda^{(\ell)} \|\bxi^{(\ell)}_{mn}\|_2 \Big] \nonumber \\
& \hspace{-.0cm} + \sum_{\ell = 1}^L \sum_{j > \ell}^L \Big[ C_{\cI(\ell,j)}(\bchi_{\cI(\ell,j)}) + \lambda_{\cI(\ell,j)} \|\bchi_{\cI(\ell,j)}\|_2 \Big] \hspace{-.2cm}  \label{Lagrangian} \\
& \hspace{-.0cm} + \sum_{\ell = 1}^L \sum_{j > \ell}^L \Big[\bmu_{\cI(\ell,j)}^\cT(\bz_{\cI(\ell,j)} - \bchi_{\cI(\ell,j)} )  + \frac{\kappa}{2} \|\bz_{\cI(\ell,j)} - \bchi_{\cI(\ell,j)}\|_2^2 \Big] \nonumber \\
& \hspace{-.0cm} + \sum_{\ell = 1}^L \hspace{-.1cm}  \sum_{j \in \cN^{(\ell)}} \hspace{-.2cm} \Big[\bgamma^{(\ell) \, \cT}_{\cI(\ell,j)}(\bxi_{\cI(\ell,j)}^{(\ell)} - \bz_{\cI(\ell,j)}) + \frac{\kappa}{2} \|\bxi_{\cI(\ell,j)}^{(\ell)} - \bz_{\cI(\ell,j)}\|_2^2 \Big]  \nonumber
\end{align}
where $\kappa > 0$ is a given constant~\cite[Sec.~3.4]{BeT89}, and $\bd :=  \{\bgamma^{(\ell)}_{\cI(\ell,j)}\}, \{\bmu_{\cI(\ell,j)}\}$ for notation brevity. 
ADMM amounts to iteratively performing the following steps ($i$ denotes the iteration index):

\noindent \textbf{[S1a]} For $\ell = 1,\ldots,L$, update $\bxi^{(\ell)}(i+1), \bs_{G}^{(\ell)}(i+1)$ as
\begin{align}
& \{ \bxi^{(\ell)}(i+1), \bs_{G}^{(\ell)}(i+1) \} =  \\
&   \arg \min_{ (\bxi^{(\ell)}, \bs_{G}^{(\ell)}) \in \cF^{(\ell)}} \cL(\bxi^{(\ell)}, \bs_{G}^{(\ell)}, \{\bchi_{\cI(\ell,j)}(i)\}, \{\bz_{\cI(\ell,j)}(i)\}, \bd(i)) \, .\nonumber 
\end{align}

\noindent \textbf{[S1b]} Per line connecting areas $\ell$ and $j$, compute
\begin{align}
&  \bchi_{\cI(\ell,j)}(i+1) = \nonumber \\
&   \arg \min_{\bchi_{\cI(\ell,j)}} \cL(\{\bxi^{(\ell)}(i+1), \bs_{G}^{(\ell)}(i+1)\}, \bchi_{\cI(\ell,j)}, \bz_{\cI(\ell,j)}(i), \bd(i)) \nonumber \\
& \textrm{s.t.~~} \quad   \bchi_{\cI(\ell,j)} \bbmM_{\cI(\ell,j)}^{\phi} \bchi_{\cI(\ell,j)} \leq (I_{\cI(\ell,j)}^{\textrm{max}})^2 , \,\, \forall \, \phi \, .
\end{align}

\noindent \textbf{[S2]} Update auxiliary variables $\{\bz_{\cI(\ell,j)}(i+1)\}$ as
\begin{align}
& \{\bz_{\cI(\ell,j)}(i+1)\} = \arg \min_{\{\bz_{\cI(\ell,j)}\}} \cL(\{\bxi^{(\ell)}(i+1), \bs_{G}^{(\ell)}(i+1)\}, \nonumber \\
& \hspace{3cm}  \{\bchi_{\cI(\ell,j)}(i+1)\}, \{\bz_{\cI(\ell,j)}\}, \bd(i) ) \, .
\end{align}

\noindent \textbf{[S3]} Update $\{\bd(i+1)\}$ using
\small
\begin{align}
& \hspace{-.4cm} \bgamma^{(\ell)}_{\cI(\ell,j)}(i+1)  = \bgamma^{(\ell)}_{\cI(\ell,j)}(i) + \kappa \big( \bxi_{\cI(\ell,j)}^{(\ell)}(i+1) - \bz_{\cI(\ell,j)}(i+1) \big)  \hspace{-.1cm} \\
& \hspace{-.4cm}  \bmu_{\cI(\ell,j)}(i+1)  = \bmu_{\cI(\ell,j)}(i) + \kappa \big( \bz_{\cI(\ell,j)}(i+1) - \bchi_{\cI(\ell,j)}(i+1)  \big) . \hspace{-.2cm}
\end{align}
 \normalsize

\noindent Step [S2] decouples into $|\cI|$ quadratic programs (one per line $\cI(\ell,j)$), each solvable in closed form. What is more, a closer look at [S2]--[S3] reveals that the dual variables satisfy the condition $\bgamma^{(\ell)}_{\cI(\ell,j)}(i) + \bgamma^{(j)}_{\cI(\ell,j)}(i) - \bmu_{\cI(\ell,j)}(i) = \mathbf{0}$ per neighboring areas $\ell$ and $j$ and iteration $i \geq 1$, whenever they are initialized at $\bgamma^{(\ell)}_{\cI(\ell,j)}(0) = \bgamma^{(j)}_{\cI(\ell,j)}(0) = \bmu_{\cI(\ell,j)}(0) = \mathbf{0}$. Using this distinct feature of the dual iterates, and leveraging the decomposability of the Lagrangian, steps [S1]--[S3] can be conveniently simplified as shown next. 

\vspace{.2cm}

\noindent \textbf{[S1a$^\prime$]} Per area $\ell = 1,\ldots,L$, solve
\begin{align}
& \{ \bxi^{(\ell)}(i+1), \bs_{G}^{(\ell)}(i+1) \} = \nonumber \\
&   \arg \min_{\bxi^{(\ell)}, \bs_{G}^{(\ell)}} \Big\{ C^{(\ell)}(\bxi^{(\ell)}, \{\bs_{G}^{(\ell)}\}) + \hspace{-.3cm} \sum_{(m,n)\in\cE_R^{(\ell)}} \hspace{-.3cm} \lambda^{(\ell)} \|\bxi^{(\ell)}_{mn}\|_2 \nonumber \\
& \hspace{.3cm} + \sum_{j \in \cN^{(\ell)}} \Big[  (\bgamma^{(\ell)}_{\cI(\ell,j)}(i))^\cT \bxi^{(\ell)}_{\cI(\ell,j)} + \frac{\kappa}{2} \|\bxi^{(\ell)}_{\cI(\ell,j)}\|_2^2  \nonumber    \\
& \hspace{1.0cm}  - \frac{\kappa}{3} (\bxi^{(\ell)}_{\cI(\ell,j)})^\cT \big(\bxi^{(\ell)}_{\cI(\ell,j)}(i) + \bxi^{(j)}_{\cI(\ell,j)}(i) + \bchi_{\cI(\ell,j)}(i)    \big) \Big] \Big\}\nonumber \\
& \textrm{s.t.~~} \quad  (\bxi^{(\ell)}, \bs_{G}^{(\ell)}) \in \cF^{(\ell)} \,  . \label{stepS1aprime}
\end{align}

\noindent \textbf{[S1b$^\prime$]} Per line $\cI(\ell,j)$, solve the constrained MSTO problem
\begin{align}
&  \bchi_{\cI(\ell,j)}(i+1) =  \arg \min_{\bchi} \Big\{ C_{\cI(\ell,j)}(\bchi) + \lambda_{\cI(\ell,j)} \|\bchi\|_2  + \frac{\kappa}{2} \|\bchi\|_2^2 \nonumber \\
& \hspace{-.0cm} - \bchi^\cT \Big[ \bmu_{\cI(\ell,j)}(i)  + \frac{\kappa}{3} \big( \bxi^{(\ell)}_{\cI(\ell,j)}(i) +  \bxi^{(j)}_{\cI(\ell,j)}(i) + \bchi_{\cI(\ell,j)}(i) \big) \Big]  \Big\} \nonumber \\
& \textrm{s.t.~~} \quad \quad   \bchi_{\cI(\ell,j)} \bbmM_{\cI(\ell,j)}^{\phi} \bchi_{\cI(\ell,j)} \leq (I_{\cI(\ell,j)}^{\textrm{max}})^2 , \,\, \forall \, \phi \, . \label{stepS1bprime}
\end{align}

\noindent \textbf{[S2$^\prime$]} Update dual variables
\begin{align}
& \hspace{-.2cm} \bgamma^{(\ell)}_{\cI(\ell,j)}(i+1)  = \bgamma^{(\ell)}_{\cI(\ell,j)}(i) \nonumber \\ 
& \hspace{-.3cm} \ + \frac{\kappa}{3} \Big(2 \bxi_{\cI(\ell,j)}^{(\ell)}(i+1) - \bxi_{\cI(\ell,j)}^{(j)}(i+1) - \bchi_{\cI(\ell,j)}(i+1) \Big) 
\label{stepS2aprime} \\
& \hspace{-.2cm} \bmu_{\cI(\ell,j)}(i+1)  = \bmu_{\cI(\ell,j)}(i) \nonumber \\ 
& \hspace{-.3cm} \ + \frac{\kappa}{3} \Big(\bxi_{\cI(\ell,j)}^{(\ell)}(i+1) + \bxi_{\cI(\ell,j)}^{(j)}(i+1) - 2 \bchi_{\cI(\ell,j)}(i+1) \Big) .  \hspace{-.1cm} \label{stepS2bprime}
\end{align}

Convergence to the solution of the centralized problem (MR3) is formalized next; see also~\cite[Sec.~3.4]{BeT89}.    

%%%%%%%%%%%%%%%%%%%%%%%%%%%%%%%%%%%%%%%
\begin{algorithm}[t]
\label{alg:decentralized}
\caption{ADMM-based distributed reconfiguration} \small{
\begin{algorithmic}

\STATE Set $\bgamma^{(\ell)}_{\cI(\ell,j)}(0) = \bgamma^{(j)}_{\cI(\ell,j)}(0) = \bmu_{\cI(\ell,j)}(0) = \mathbf{0}$ for all $\ell, j$.

\FOR {$i = 1,2,\ldots$ (repeat until convergence)} 

\STATE 1. [LAC $\ell$]: compute $\bxi^{(\ell)}(i+1)$ and $\bs_{G}^{(\ell)}(i+1)$ via~\eqref{stepS1aprime}. \\

\hspace{.25cm}  [MGM]: compute $\{\bchi_{\cI(\ell,j)}(i+1)\}$ via~\eqref{stepS1bprime}. \\

\STATE 2. [LAC $\ell$]: send $\{\bxi_{\cI(\ell,j)}^{(\ell)}(i+1), j \in \cN^{(\ell)} \}$ to MGM. 

\hspace{.25cm}  [MGM]: receive $\{\bxi_{\cI(\ell,j)}^{(\ell)}(i+1), j \in \cN^{(\ell)} \}$ from LAC $\ell$. 

\hspace{1.45cm}  Repeat for all $\ell = 1, \ldots, L$. 

\hspace{.25cm}  [MGM]: send $\bchi_{\cI(\ell,j)}(i+1)$ and $\{\bxi_{\cI(\ell,j)}^{(j)}(i+1), \forall \, j \in \cN^{(\ell)}\}$. 

\hspace{1.45cm} Repeat for all $\ell = 1, \ldots, L$.

\STATE 3. [LAC $\ell$]: update dual variables $\{\bgamma^{(\ell)}_{\cI(\ell,j)}(i+1), j \in \cN^{(\ell)}\}$.

\hspace{.25cm}  [MGM]: update dual variables $\{\bmu_{\cI(\ell,j)}(i+1), \cI(\ell,j) \in \cI\}$.

\ENDFOR
\end{algorithmic}}
\end{algorithm}
%%%%%%%%%%%%%%%%%%%%%%%%%%%%%%%%%%%%%%

\vspace{.2cm}

\begin{proposition}
Suppose that the dual variables are initialized at $\bgamma^{(\ell)}_{\cI(\ell,j)}(0) = \bgamma^{(j)}_{\cI(\ell,j)}(0) = \bmu_{\cI(\ell,j)}(0) = \mathbf{0}$. Then, for any $\kappa > 0$ the iterates $\{\bxi^{(\ell)}(i), \bs_{G}^{(\ell)}(i), \bchi_{\cI(\ell,j)}(i), \bd(i)\}$ obtained from [S1$^\prime$]--[S2$^\prime$] are convergent, and $\lim_{i \rightarrow + \infty} \{\bxi^{(\ell)}(i), \bs_{G}^{(\ell)}(i), \ell = 1,\ldots,L\} = \{\bxi_\mathrm{opt}, \bsigma_{G_n,\mathrm{opt}}^\phi\}$, with $\{\bxi_\mathrm{opt}, \bsigma_{G_n,\mathrm{opt}}^\phi\}$ the optimal solution of (MR3).  \hfill $\Box$
\end{proposition}

\vspace{.2cm}

The resulting decentralized algorithm entails a two-way communication between the MGM and the LACs, and it is tabulated as Algorithm 1. Per iteration $i$, each LAC updates $\bxi^{(\ell)}(i+1)$ and $\bs_{G}^{(\ell)}(i+1)$ via [S1a$^\prime$], and subsequently sends the $|\cN^{(\ell)}|$ real-valued vectors $\{\bxi_{\cI(\ell,j)}^{(\ell)}(i+1), j \in \cN^{(\ell)} \}$ to the MGM. After performing step [S1b$^\prime$], the MGM receives $\{\bxi_{\cI(\ell,j)}^{(\ell)}(i+1), j \in \cN^{(\ell)} \}$ from each LAC, and then it sends $\bchi_{\cI(\ell,j)}(i+1)$ along with $\{\bxi_{\cI(\ell,j)}^{(j)}(i+1), \forall \, j \in \cN^{(\ell)}\}$ to the LAC $\ell$, for $\ell = 1, \ldots, L$. Finally, MGM and LACs update the dual variables via [S2$^\prime$]. The overall distributed procedure is tabulated as Algorithm 1.

\vspace{.1cm}

\emph{Remark 3}. The premise of Algorithm 1 is that samples of the solar, wind, and load forecasting errors are available, and they are employed to find the quantities on the right hand side of~\eqref{approx_solp2}. Errors in the forecasts of solar irradiance and wind speed may be correlated across space~\cite{Bacher2009,Lorenz09,Boone05}, especially for geographically close RES facilities. Load forecasting errors can be roughly approximated as spatially uncorrelated~\cite{Chertkov13}. There are three viable setups where spatially-correlated samples can be  generated, depending on the role of  the MGM:

\emph{s1)} Statistics of the prediction errors are readily available when forecasts are carried out at the MGM for the entire microgrid;  the MGM performs Monte Carlo sampling, and subsequently disseminates a $2 \times 1$ real-valued vector per phase and node to be used in~\eqref{approx_solp2}.

\emph{s2)}  Each LAC performs the forecasting for its own area, and notifies the MGM about the prediction errors. The empirical joint distribution of  the prediction errors is obtained at the MGM, which computes the quantities in~\eqref{approx_solp2}.   

\emph{s3)} Forecasts are performed at the LACs for their own areas, and synthetic spatial correlation models are used to draw the samples; see e.g.,~\cite{Lorenz09,Boone05}.

In specific setups, solar, wind, and load forecasting errors can be spatially uncorrelated as explained in~\cite{Chertkov13}.

%%%%%%%%%%%%%%%%%%%%%%%%%%%%%%%%%%%%%%%%%%%%%%
\section{Numerical Experiments}
\label{sec:Simulation}
%%%%%%%%%%%%%%%%%%%%%%%%%%%%%%%%%%%%%%%%%%%%%%

\begin{figure}[t]
\begin{center}
\includegraphics[width=0.50\textwidth]{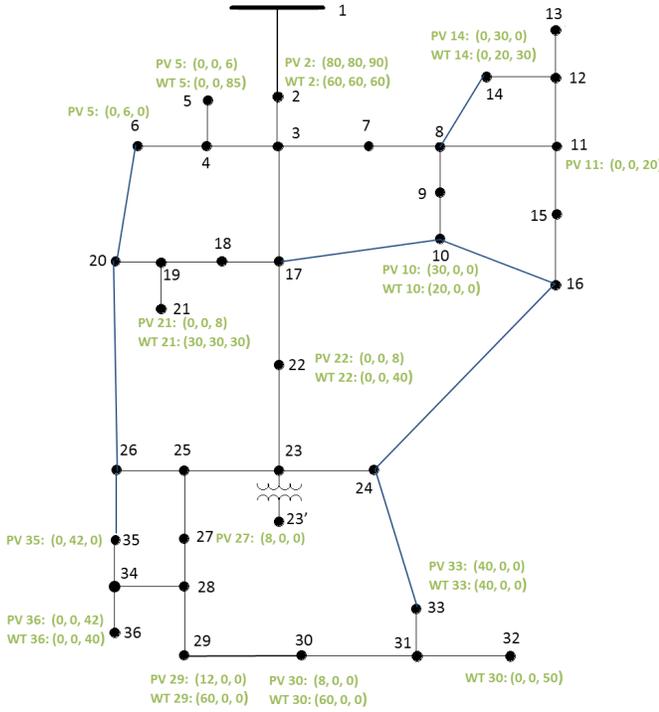}
\vspace{-.2cm}
\caption{Modified IEEE 37-bus test feeder. }
\label{fig:F_feeder37}
\vspace{-.2cm}
\end{center}
\end{figure}

The effectiveness of the proposed scheme is showcased in this section using a modified version of the IEEE 37-node test feeder~\cite{testfeeder}. As shown in Fig.~\ref{fig:F_feeder37},  
eight three-phase lines are added to the original radial scheme, and DERs are placed throughout the network. The parameters of the additional lines are listed in Table~\ref{tab:lines}, where the line matrices corresponding to the configuration indexes 723 and 724 can be found in~\cite{testfeeder}. Further, the 17 branches $\cE_{R} = \{(1,2), (3,4), (6,20), (7,8), (8,9), (8,14), (15,16), (16,24),$ $(10,16), (10,17), (17,18), (20,26), (23,24),(23,25),(24,33),$ $(29,30), (26,35)\}$ feature sectionalizing switches. Controllable DG units are located at nodes $\{10, 12, 16, 19, 24, 28, 32\}$, they operate at unity power factor, and they can supply a maximum power of $50$ kW per phase. 
PV systems and small wind turbines (WTs) operate at unity power factor, and their generation capacity per phase $(a,b,c)$ and node (in kW-peak) is reported in Fig.~\ref{fig:F_feeder37}. Finally, the impedance matrices for the original lines and the loads are the ones specified in~\cite{testfeeder}.

The package \texttt{CVX}\footnote{[Online] Available: \texttt{http://cvxr.com/cvx/}} is used to solve the reconfiguration problem in \texttt{MATLAB}.  The average computational time required by the interior-point solver of \texttt{CVX} was $0.8$ seconds on a machine with Intel Core i7-2600 CPU @ 3.40GHz. 

%%%%%%%%%%%%%%%%%%%%
%%%%%%%%%%%%%%%%%%%%

%%%%%%%%%%%%%%%%%%%%
\begin{table}[t]
\caption{Additional lines in the modified IEEE 37-node feeder}
\begin{center}
\begin{tabular}{|c|cc||c|cc}
Line & Conf. & Length (ft) & Line & Conf. & Length (ft) \\
\hline
(8,14) & 723 & 1144 & (16,24) & 724 & 1580\\
(6,20) & 724 & 1320  & (10,17) & 724 & 1137  \\
(10,16) & 724 & 847 & (24,33) & 724 & 1315 \\
(20,26) & 724 & 815 & (26,35) & 724 & 377 \\
\end{tabular}
\vspace{-.8cm}
\end{center}
\label{tab:lines}
\end{table}%
%%%%%%%%%%%%%%%%%%%%

To account for forecasting errors, the actual power supplied by RESs is modeled as $P_{E_n}^{\phi} = \bar{P}_{E_n}^{\phi} + \Delta_{E_n}^{\phi}$, with $\bar{P}_{E_n}^{\phi}$ the (known) forecasted value and $\Delta_{E_n}^{\phi}$ the (random) forecasting error. A zero-mean truncated Gaussian distribution is adopted for $\Delta_{E_n}^{\phi}$, with truncation at the $0.13$th and $99.87$th percentiles; see e.g.,~\cite{Bacher2009,Tsikalakis06}. Random variables $\{\Delta_{E_n}^{\phi}\}$ are correlated across nodes, and their correlation matrix is obtained using an exponentially decreasing function of the distance between nodes as specified in~\cite{Lorenz09} and~\cite[Ch.~9]{Boone05} for PV systems and WTs, respectively. 
Load forecasting errors are modeled as $S_{L_n}^{\phi} = \bar{S}_{L_n}^{\phi} + (\Delta_{L_n,P}^{\phi} + j \Delta_{L_n,Q}^{\phi})$, with $\bar{S}_{L_n}^{\phi}$ denoting the forecasted value, and $\Delta_{L_n,P}^{\phi}, \Delta_{L_n,Q}^{\phi}$ capturing errors in the prediction of the active and reactive loads, respectively. Variables $\{\Delta_{L_n,P}^{\phi}\}, \{\Delta_{L_n,Q}^{\phi}\}$  
are Gaussian distributed~\cite{Hodge2013}, zero-mean, uncorrelated, and truncated at the $0.13$th and $99.87$th percentiles. The distribution of the approximation errors $\{\bepsilon_{\iota}^{\phi}\}$ was evaluated via extensive simulations, by comparing the injected currents obtained from (MR3) without error compensation, with the ones obtained via OPF~\cite{Dallanese-TSG13}.   

%%%%%%%%%%%%%%%%%%%%%
\begin{figure}
	\centering
  \subfigure[Setup 1: the standard deviation of the solar power prediction error amounts to $5 \%$ of the forecasted value; the error on the wind power is on the order of $20 \%$ of the forecasted value; and, the standard deviation of the load forecasting error is in the interval $4-6 \%$ of $\bar{S}_{L_n}^\phi$]{\includegraphics[width=0.48\textwidth]{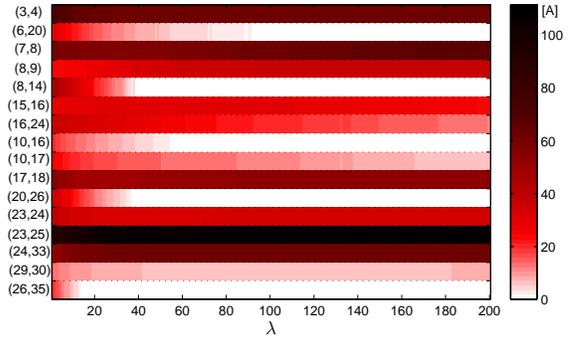}}
  \subfigure[Setup 2: the standard deviations of the solar power, wind power, and load are set to $0.05 \%$, $0.2 \%$, and $0.4-0.5 \%$, respectively.]{\includegraphics[width=0.48\textwidth]{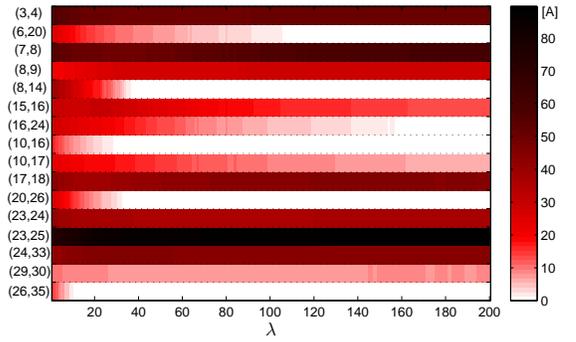}}
  \caption{Sum of current magnitudes $\sum_{\phi \in \cP_{mn}}|I_{mn}^\phi|$ on lines $\cE_R$.}
  \label{fig:currents}
  \vspace{-.5cm}
\end{figure}
%%%%%%%%%%%%%%%%%%%%

Let $I_{mn}^{\textrm{max}} = 300$ A for the conductors of line $(1,2)$; $150$ A for lines $(2,3), (3,17)$; and, $100$ A for all the remaining branches. Consider the group-sparsity regularization [cf.~\eqref{Glasso_currents}]
$g(\bxi_R) :=  \sum_{(m,n) \in \cE_{R}} \, \lambda \omega_{mn} \| \bxi_{mn} \|_2$ where $\omega_{mn} = 1$ for the lines in the original feeder scheme~\cite{testfeeder}, and  $\omega_{mn} = 1.5$ for the $8$ additional lines in Table~\ref{tab:lines}. This way, utilization of the original lines in encouraged. Suppose that for a given number of utilized lines, the goal is to minimize the net microgrid operational cost; that is, $C(\cV) = \sum_{\phi} \Re\{V_1^\phi (I_1^\phi)^{*}\} + \hspace{-.2cm} 0.5 \sum_{n, \phi } P_{G,n}^{\phi} + \sum_{(m,n) \in \cE} \Re\{ \bi_{mn}^{\cT} \bZ_{mn} \bi_{mn}\}$. Suppose that the forecasted solar power amounts to $90 \%$ of the kW-peak, while the WTs are operating at $70 \%$ of their maximum capacity. The threshold for the LOL probability is set to $\rho = 0.01$ ($1 \%$), and the parameter $\beta$ in~\eqref{sample_size_MR3} is $0.05$. 

Two setups are considered: 

\emph{Setup 1}: the standard deviation of the solar power prediction error amounts to $5 \%$ of the forecasted value~\cite{Bacher2009}; the error on the wind power is on the order of $20 \%$ of the forecasted value~\cite{Tsikalakis06}; and, the standard deviation of the load forecasting error is in the interval $[4,6] \%$ of $\bar{S}_{L_n}^\phi$~\cite{Hodge2013}; 

\emph{Setup 2}: the standard deviations are set to $0.05 \%$, $0.2 \%$, and $[0.4, 0.5] \%$ to resemble a markedly higher prediction accuracy~\cite{Bacher2009,Lorenz09,Tsikalakis06}. 

For given load conditions (or, minimum load requirements~\eqref{approx_solp2}), optimality and complexity of the sparsity-based reconfiguration scheme were already discussed in~\cite{Dallanese-TPD13}. Here, the objective is to \emph{i)} investigate the effects of load and renewable generation uncertainties on the optimal topology, and \emph{ii)} assess convergence of the proposed decentralized reconfiguration protocol. 

Fig.~\ref{fig:currents} depicts the sum of the current magnitudes $\sum_{\phi \in \cP_{mn}}|I_{mn}^\phi|$ on lines equipped with switches, for different values of the tuning parameter $\lambda$. 
The current magnitude is color-coded, where white represents zero current; this means that the switches are open, and the distribution line is not utilized. First, notice that the number of open switches increases as $\lambda$ increases, thus further corroborating the results in~\cite{Dallanese-TPD13}. By varying $\lambda$, the MGM can obtain either meshed topologies (low values of $\lambda$), weakly-meshed (high values of $\lambda$), or even radial by simply adjusting $\lambda$. 
Notice however that, in the simulated setups, it is  impossible to find a tree topology for which loss-of-load probabilities and Ampacity limits are satisfied. Take for example $\lambda = 200$: 5 switches are open in Fig.~\ref{fig:currents}(a), and 6 in Fig.~\ref{fig:currents}(b),  and thus weakly meshed topologies are obtained. Numerical experiments reveal that if one ``opens'' the switch on line $(10,17)$ in the first setup, then (MR3) is infeasible for many realizations of the forecasting errors. 
Comparing Figs.~\ref{fig:currents}(a) and (b), it is observed that for a fixed value of $\lambda$, the number of lines utilized in the first setup is typically higher that in the second. Indeed, when the forecasting error is high, it is prudent to utilize a higher number of lines to avoid exceeding Ampacity limits if the actual RES generation and load demand deviate from the forecasted values. 

Finally, if the wanted topology is weakly-meshed with a pre-specified number of lines to utilize, the distribution system operator can quickly gauge the optimal topology from printouts like Fig.~\ref{fig:currents}, for given values of the forecasting error standard deviation. 

This example highlights the merits of the proposed risk-constrained reconfiguration approach. Specifically, the obtained topology is: \emph{i)} optimal according to the regularized optimization criterion $C(\cV)$, rather than being a result of line selection heuristics~\cite{Khodr09,Celli02,Baran89,Schmidt95}, which are computationally heavy and may identify sub-optimal configurations; and, \emph{ii)} it guarantees feasible power flow solutions for the majority of the (unknown) RES generation and load realizations.              

\begin{figure}[t]
\begin{center}
\includegraphics[width=0.50\textwidth]{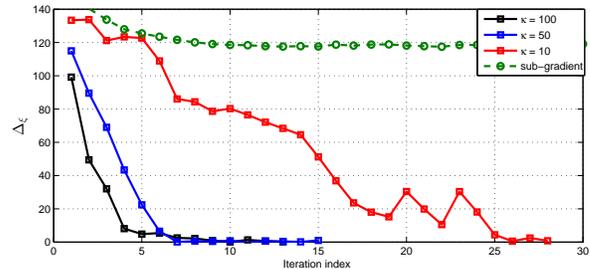}
\caption{Convergence of the distributed scheme. }
\label{fig:F_convergence}
\vspace{-.5cm}
\end{center}
\end{figure}

Finally, convergence of the ADMM-based decentralized algorithm is exemplified in Fig.~\ref{fig:F_convergence}, where $3$ areas within the microgrid are managed autonomously, while the rest of the network is controlled by the MGM. Specifically, the three areas are formed by the subsets of nodes $\cA^{(1)} = \{11, 12, 13, 14, 15\}$, $\cA^{(2)} = \{18, 19, 20, 21\}$, and $\cA^{(3)} = \{4, 5, 6\}$. As a representative example, the trajectories corresponding to  $\Delta_{\xi}(i) := \|\bxi_{1,2}^{(1)}(i) - \bxi_{1,2}^{(2)}(i)\|$, with $\bxi_{1,2}^{(\ell)}(i)$ stacking real and imaginary parts of the currents on lines $(8,14), (8,11)$, and $(15,16)$ per iteration $i$, are reported for different values of the ADMM parameter $\kappa$. They are also 
compared with the ones obtained by using the sub-gradient ascent-based distributed algorithm, with constant stepsize of $0.1$~\cite{Nedic09}. The proposed distributed solver exhibits markedly faster convergence than the one based on the sub-gradient. It can be seen that the discrepancies between $\bxi_{1,2}^{(1)}(i)$ and $\bxi_{1,2}^{(2)}(i)$ vanish more rapidly as $\kappa$ increases, which is in line with the discussion in e.g.,~\cite{BoydADMoM,ErsegheADMM}, where it is pointed out that relatively large values of $\kappa$ place a large penalty on violations of primal feasibility. In general,  ADMM convergence rate depends not only on $\kappa$, but also on the strong convexity constant of the cost functions, the Lipschitz constants of their gradients, and the condition number of matrix $\bA$, used to rewrite constraints~\eqref{P3bconsensus}--\eqref{P3bconsensus2}  in compact form as $\bA \bnu = \mathbf{0}$, for an appropriately defined vector $\bnu$; see e.g.,~\cite{Deng12}. For simpler problem setups, given the condition number of  $\bA$ and the other quantities related to the cost function, it is possible to find the value of $\kappa$ that maximizes the convergence speed~\cite{Deng12}. Extending the results of~\cite{Deng12} to the present setup, along with analyzing the ties between $\kappa$ and convergence rate will be the subject of future research. Notice, however, that computing Lipschitz constants and condition numbers requires global knowledge of electrical network and per-area optimization objectives; thus, it may not be feasible in the present setup.

%%%%%%%%%%%%%%%%%%%%
%%%%%%%%%%%%%%%%%%%%

%%%%%%%%%%%%%%%%%%%%%%%%%%%%%%%%%%%%%%%%%%%%%%
\section{Concluding Remarks}
\label{sec:conclusions}
%%%%%%%%%%%%%%%%%%%%%%%%%%%%%%%%%%%%%%%%%%%%%%

The system reconfiguration task was considered for microgrids, in the presence of renewable-based generation and load foresting errors. To cope with possible supply-demand imbalance, a novel chance-constrained optimization problem was formulated to limit the probability of LOL, while adhering to line Ampacity constraints strictly. The novel reconfiguration approach utilizes sparsity-promoting regularization terms to effect line selection, and a scenario optimization technique to 
approximate the probabilistic constraints. The upshot of the proposed formulation is that it leads to a convex program, and it entails one balance constraint per phase and node.  Finally, a novel decentralized reconfiguration scheme was developed, which entails a two-way communication between the MGM and the LACs to consent on the value of the currents flowing on the lines interconnecting the areas.

%%%%%%%%%%%%%%%%%%%%%%%%%%%%%%%%%%%%%%%%%%%%%%
\bibliographystyle{IEEEtran}
\bibliography{biblio.bib}
%%%%%%%%%%%%%%%%%%%%%%%%%%%%%%%%%%%%%%%%%%%%%%

%%%%%%%%%%%%%%%%%%%%%%%%%%%%%%%%%%%%%%%%%%%%%%
% biography
%%%%%%%%%%%%%%%%%%%%%%%%%%%%%%%%%%%%%%%%%%%%%%

\begin{biography}[{\includegraphics[width=1in,height=1.25in,clip,keepaspectratio]{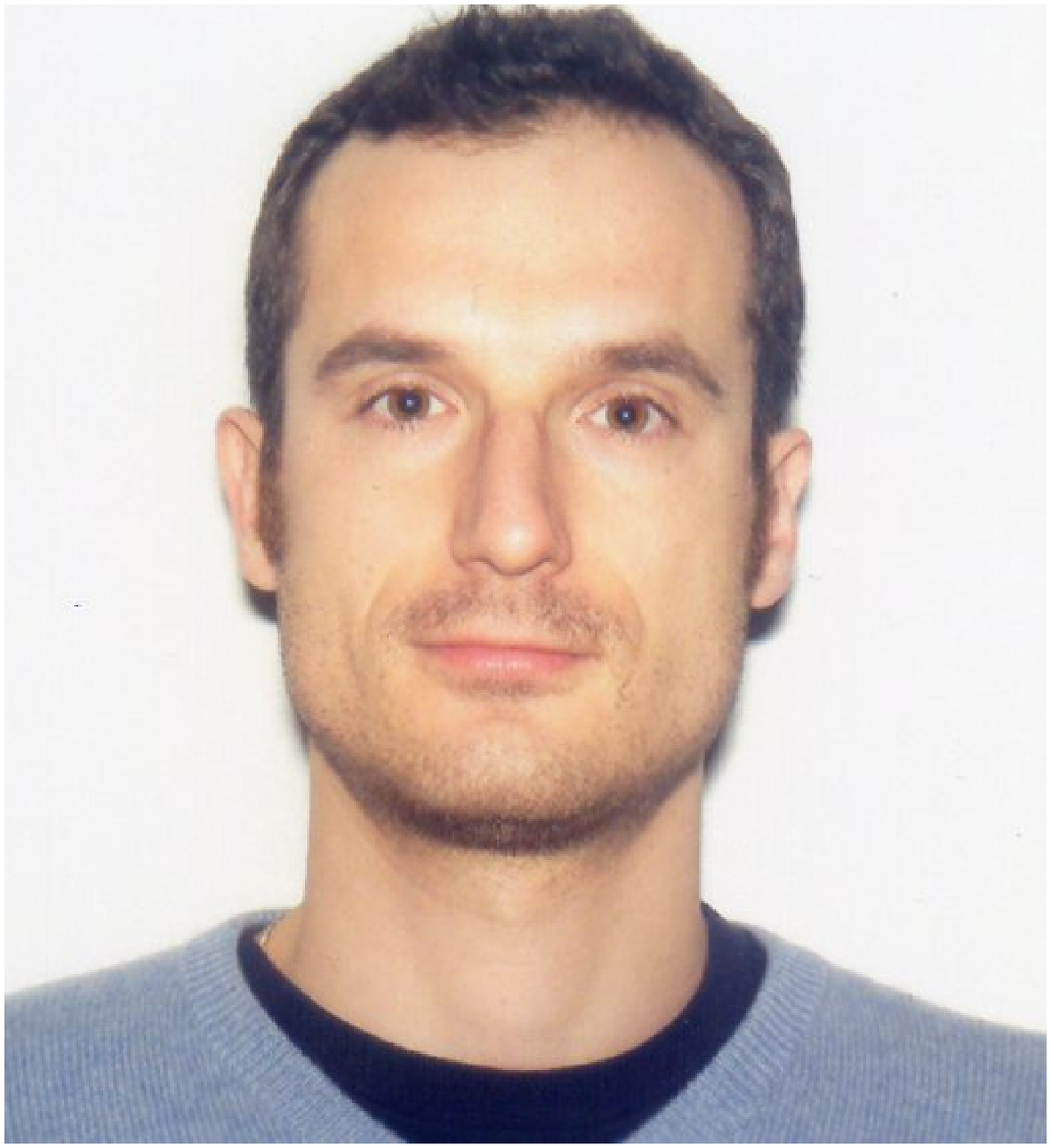}}]
{Emiliano Dall'Anese (S'08-M'11)} received the Laurea Triennale (B.Sc Degree) and the Laurea Specialistica (M.Sc Degree) from the University of Padova, Italy, in 2005 and 2007, respectively, and the Ph.D. in Information Engineering from the Department of Information Engineering, University of Padova, Italy, in 2011. From January 2009 to September 2010, he was a visiting scholar at the Department of Electrical and Computer Engineering, University of Minnesota, USA. Since January 2011, he has been a Postdoctoral Associate at the Department of Electrical and Computer Engineering and Digital Technology Center, University of Minnesota, USA. 

His research interests lie in the areas of Power Systems, Signal Processing, and Communications.
Current research focuses on the analysis and application of  Signal Processing, Optimization, and
Machine Learning tools and methods to energy management in future power systems and grid informatics.
\end{biography}

\begin{biography}[{\includegraphics[width=1in,height=1.25in,clip,keepaspectratio]{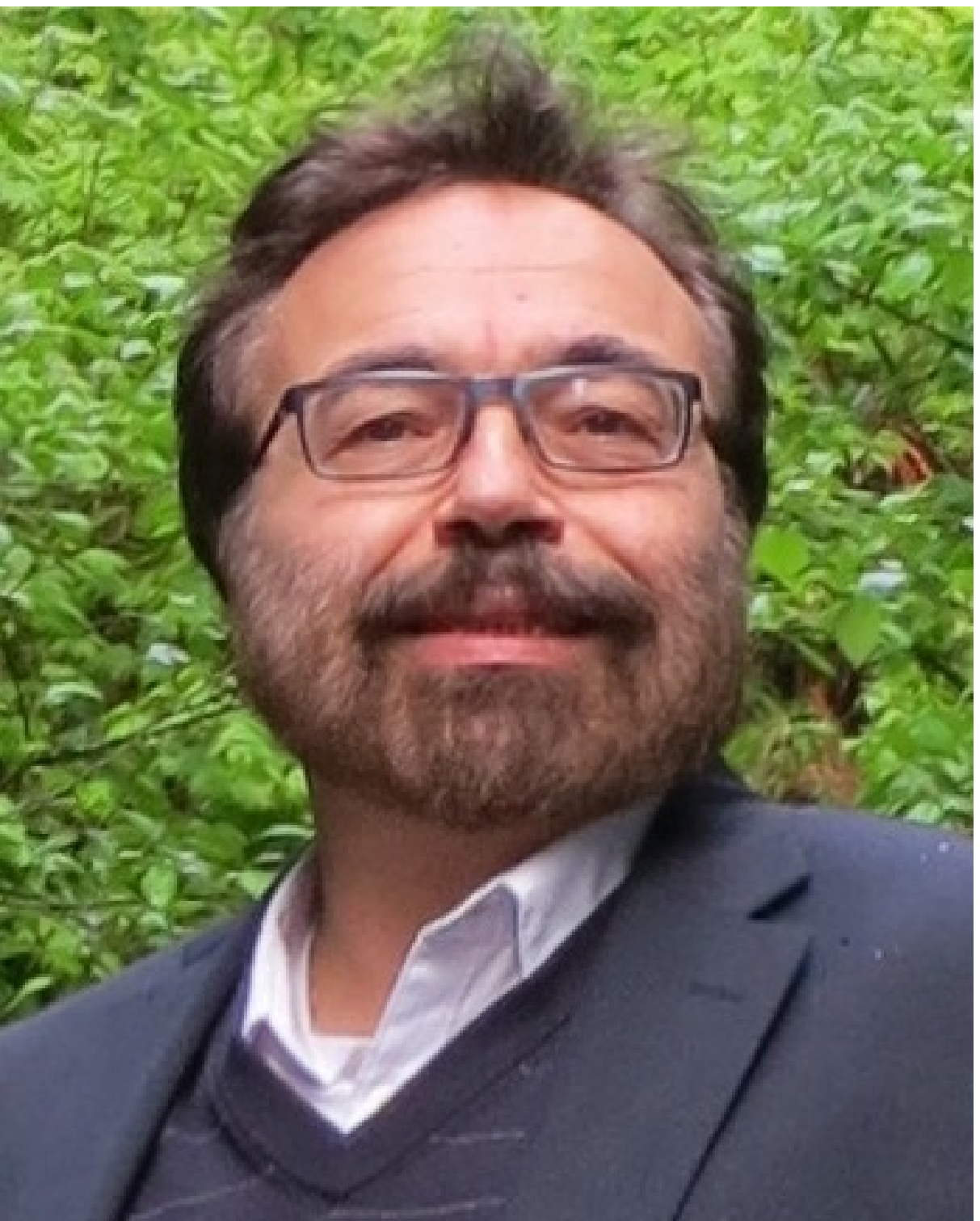}}]
{Georgios B.  Giannakis (Fellow'97)} received his Diploma in Electrical Engr. from the Ntl. Tech. Univ. of Athens, Greece, 1981. From 1982 to 1986 he was with the Univ. of Southern California (USC), where he received his MSc. in Electrical Engineering, 1983, MSc. in Mathematics, 1986, and Ph.D. in Electrical Engr., 1986. Since 1999 he has been a professor with the Univ. of Minnesota, where he now holds an ADC Chair in Wireless Telecommunications in the ECE Department, and serves as director of the Digital Technology Center.

His general interests span the areas of communications, networking and statistical signal processing - subjects on which he has published more than 360 journal papers, 600 conference papers, 21 book chapters, two edited books and two research monographs (h-index 107). Current research focuses on sparsity and big data analytics, wireless cognitive radios, mobile ad hoc networks, renewable energy, power grid, gene-regulatory, and social networks. He is the (co-) inventor of 22 patents issued, and the (co-) recipient of 8 best paper awards from the IEEE Signal Processing (SP) and Communications Societies, including the G. Marconi Prize Paper Award in Wireless Communications. He also received Technical Achievement Awards from the SP Society (2000), from EURASIP (2005), a Young Faculty Teaching Award, and the G. W. Taylor Award for Distinguished Research from the University of Minnesota. He is a Fellow of EURASIP, and has served the IEEE in a number of posts, including that of a Distinguished Lecturer for the IEEE-SP Society.
\end{biography}

\end{document}